\documentclass[11pt,a4paper,fleqn]{article}
\usepackage{a4wide,amsfonts,amsmath,latexsym,amssymb,euscript,graphicx,units,mathrsfs}

\usepackage{graphicx}
\usepackage{color}
\usepackage{amssymb}
\usepackage{amssymb}
\usepackage[T1]{fontenc}
\usepackage{latexsym}
\usepackage{xypic}
\usepackage{eufrak}
\usepackage{euscript}
\usepackage{amsfonts,amsmath}
\usepackage{verbatim}
\usepackage{fancyhdr}
\usepackage[english]{babel}
\usepackage{mathrsfs}
\usepackage{units}
\usepackage{stmaryrd}

\newtheorem{prop}{Proposition}[section]
\newtheorem{cor}[prop]{Corollary}
\newtheorem{lemme}[prop]{Lemma}
\newtheorem{rem}[prop]{Remark}

\newtheorem{thm}[prop]{Theorem}

\newtheorem{example}[prop]{Example}

\renewcommand{\geq}{\geqslant}
\def\leq{\leqslant}

\newcommand{\R}{\mathbb{R}}

\def\cal{\mathcal}
\def\1{{\mathbf{1}}}

\def\1{{\mathbf{1}}}
\def\0.5{{\frac{1}{2}}}

\def\H{\EuFrak H}

\newcommand{\fin}
{ \vspace{-0.6cm}
\begin{flushright}
\mbox{$\Box$}
\end{flushright}
\noindent }

\newcommand{\qed}{\nopagebreak\hspace*{\fill}
{\vrule width6pt height6ptdepth0pt}\par}

\begin{document}

\begin{center}
{\large{\bf Stein's method and normal approximation of Poisson functionals}}\\~\\
by G. Peccati \footnote{ Laboratoire de Statistique Th\'eorique et
Appliqu\'ee, Universit\'e Pierre et Marie Curie, 175 rue du Chevaleret, 75013 Paris, France, E-mail:
{\tt giovanni.peccati@gmail.com}}, J.L. Solé \footnote{Departament de Matem\`{a}tiques, Facultat de Ci\`{e}ncies, Universitat Aut\'{o}noma de Barcelona, 08193 Bellaterra
(Barcelona), Spain; E-mail: {\tt jllsole@mat.uab.es}}, M.S. Taqqu \footnote{
Boston University, Departement of Mathematics, 111 Cummington Road, Boston
(MA), USA. E-mail: \texttt{murad@math.bu.edu}. } and F. Utzet \footnote{Departament de Matem\`{a}tiques, Facultat de Ci\`{e}ncies, Universitat Aut\'{o}noma de Barcelona, 08193 Bellaterra
(Barcelona), Spain; E-mail: {\tt utzet@mat.uab.es}}\\
%\small{September 25, 2007}\\
\small{{\it This version}: July 31, 2008}
\end{center}

{\small \noindent {\bf Abstract:} We combine Stein's method with a
version of Malliavin calculus on the Poisson space. As a result, we
obtain explicit Berry-Ess\'{e}en bounds in Central Limit Theorems
(CLTs) involving multiple Wiener-It\^{o} integrals with respect to
a general Poisson measure. We provide several applications to
CLTs related to Ornstein-Uhlenbeck Lévy processes.}
\\



\section{Introduction}
In a recent series of papers, Nourdin and Peccati \cite{NouPe07,
NouPe08}, and Nourdin, Peccati and R\'{e}veillac \cite{NouPeRev},
have shown that one can effectively combine Malliavin calculus on a
Gaussian space (see e.g. \cite{NuaBook}) and Stein's method (see
e.g. \cite{Chen_Shao_sur, Reinert_sur}) in order to obtain explicit
bounds for the normal and non-normal approximation of smooth
functionals of Gaussian fields.

\smallskip

The aim of the present paper is to extend the analysis initiated in
\cite{NouPe07} to the framework of the normal approximation (in the
Wasserstein distance) of regular functionals of Poisson measures
defined on abstract Borel spaces. As in the Gaussian case, the main
ingredients of our analysis are the following:
\begin{itemize}
\item[1.] A set of \textsl{Stein differential equations}, relating
the normal approximation in the Wasserstein distance to first order
differential operators.
\item[2.] A (Hilbert space-valued) \textsl{derivative operator} $D$, acting on
real-valued square-integrable random variables.
\item[3.] An \textsl{integration by parts formula}, involving the adjoint operator of $D$.
\item[4.] A ``pathwise representation'' of $D$ which, in the Poisson
case,
involves standard \textsl{difference operators}.
\end{itemize}

As a by-product of our analysis, we obtain substantial
generalizations of the Central Limit Theorems (CLTs) for functionals
of Poisson measures (for instance, for sequences of single and
double Wiener-It\^{o} integrals) recently proved in \cite{PeTaq0,
PecTaqINT, PecTaq07} (see also \cite{DBPP, PePrU} for applications
of these results to Bayesian non-parametric statistics). In
particular, one of the main results of the present paper (see
Theorem \ref{T : chaosCLT} below) is a CLT for sequences of
multiple Wiener-It\^{o} integrals of arbitrary (fixed) order with
respect to a general Poisson measure. Our conditions are expressed in terms of
\textsl{contraction operators}, and can be seen as a Poisson
counterpart to the CLTs on Wiener space proved by Nualart and
Peccati \cite{NP} and Nualart and Ortiz-Latorre \cite{NO}.
The reader is referred to Decreusefond and Savy \cite{Dec_Savy} for
other applications of Stein-type techniques and Malliavin
operators to the assessment of Rubinstein distances on configuration
spaces.

\smallskip

The remainder of the paper is organized as follows. In Section
\ref{S : DEF}, we discuss some preliminaries, involving
multiplication formulae, Malliavin operators, limit theorems and
Stein's method. In Section \ref{S : General Ineq}, we derive a
general inequality concerning the Gaussian approximation of regular
functionals of Poisson measures. Section \ref{S : MultipleWeak} is
devoted to upper bounds for the Wasserstein distance, and Section \ref{SS : CLTchaos}
to CLTs for multiple Wiener-It\^{o}
integrals of arbitrary order. Section \ref{SS : SDint} deals with
sums of a single and a double integral. In Section \ref{S : EX},
we apply our results to non-linear functionals of Ornstein-Uhlenbeck Lévy processes.

\section{Preliminaries \label{S : DEF}}
\subsection{Poisson measures} \label{SS : PoissMeas} \setcounter{equation}{0}
Throughout the paper, $\left( Z,\mathcal{Z},\mu \right) $ indicates
a measure space such that $Z$ is a Borel space and $\mu $ is a
$\sigma$-finite non-atomic Borel measure. We define the
class $\mathcal{Z}_{\mu }$ as $\mathcal{Z}_{\mu }=\left\{ B\in \mathcal{Z}%
:\mu \left( B\right) <\infty \right\} $. The symbol $\widehat{N}=\{\widehat{N%
}\left( B\right) :B\in \mathcal{Z}_{\mu }\}$ indicates a \textsl{compensated
Poisson random measure} on $\left( Z,\mathcal{Z}\right) $ with control $\mu
. $ This means that $\widehat{N}$ is a collection of random variables
defined on some probability space $\left( \Omega ,\mathcal{F},\mathbb{P}%
\right) $, indexed by the elements of $\mathcal{Z}_{\mu }$, and such that:
(i) for every $B,C\in \mathcal{Z}_{\mu }$ such that $B\cap C=\varnothing $, $%
\widehat{N}\left( B\right) $ and $\widehat{N}\left( C\right) $ are
independent, (ii) for every $B\in \mathcal{Z}_{\mu }$,%
\begin{equation*}
\widehat{N}\left( B\right) \overset{\text{law}}{=}\mathfrak{P}\left(
B\right) -\mu \left( B\right) \text{,}
\end{equation*}%
where $\mathfrak{P}\left( B\right) $ is a Poisson random variable
with parameter $\mu \left( B\right) $. Note that properties
(i)-(ii) imply, in particular, that $\widehat{N}$ is an
\textsl{independently scattered }(or \textsl{completely random})
measure (see e.g. \cite{PecTaq07}).

\begin{rem}
{\rm According e.g. to \cite[Section 1]{Picard1996}, and due to the assumptions on the space $\left( Z,\mathcal{Z},\mu
\right)$, it is always possible to define $\left( \Omega
,\mathcal{F},\mathbb{P}\right)$ in such a way that
\begin{equation}\label{NV1}
\Omega = \left \{\omega = \sum_{j=0}^n \delta_{z_j}, \,\,
n\in\mathbb{N}\cup \{\infty\}, \,\, z_j\in Z \right \}
\end{equation}
where $\delta_z$ denotes the Dirac mass at $z$, and $\widehat{N}$
is the \textsl{compensated canonical mapping} given by
\begin{equation}\label{CanMap}
\widehat{N}(B)(\omega) = \omega(B) - \mu (B), \,\,\, B \in
\mathcal{Z}_\mu.
\end{equation}
Also, in this case one can take $\mathcal{F}$ to be the
$\mathbb{P}$-completion of the $\sigma$-field generated by the
mapping (\ref{CanMap}). Note that, under these assumptions, one has
that
\begin{equation}\label{PP}
\mathbb{P}\{\omega : \omega(B)<\infty, \,\, \forall B \text{ \ s.t. \
} \mu(B)<\infty \} =1 \,\,\,\,\text{and}\,\,\,\, \mathbb{P}\{\omega
: \exists z \,\,\text{s.t.}\,\, \omega(\{z\})>1 \} =0 .
\end{equation}
}
\end{rem}

\begin{rem}{\rm From now on, and for the rest of the paper, the hypotheses
(\ref{CanMap}) and (\ref{PP}), on the structure of $\Omega$ and
$\widehat{N}$, are implicitly satisfied. In particular,
$\mathcal{F}$ is the $\mathbb{P}$-completion of the $\sigma$-field
generated by the mapping in (\ref{CanMap}). }
\end{rem}

\medskip

For every deterministic function $h\in L^{2}\left(
Z,\mathcal{Z},\mu \right) =L^{2}\left( \mu \right) $ we write
$\widehat{N}\left( h\right) =\int_{Z}h\left( z\right)
\widehat{N}\left( dz\right) $ to indicate the
Wiener-It\^{o} integral of $h$ with respect to $\widehat{N}$ (see e.g. \cite{KypBook} or \cite{Sato}%
). We recall that, for every $h\in L^{2}\left( \mu \right) $, the random variable $\widehat{N}%
\left( h\right) $ has an infinitely divisible law, with L\'{e}vy-Khinchine
exponent (see again \cite{Sato}) given by%
\begin{equation}
\psi \left( h,\lambda \right)
=\log\mathbb{E}\left[{\rm e}^{i\lambda\widehat{N}%
\left( h\right) }\right]=\int_{Z}\left[{\rm e} ^{i\lambda h\left(
z\right)}
-1-i\lambda h\left( z\right)\right] \mu ( dz), \text{ \ \ }%
\lambda \in \mathbb{R}.  \label{PLK}
\end{equation}%
Recall also the isometric relation: for every $g,h\in L^{2}\left( \mu
\right) $, $\mathbb{E[}\widehat{N}\left( g\right) \widehat{N}\left( h\right)
]=\int_{Z}h\left( z\right) g\left( z\right) \mu \left( dz\right) $.

\medskip

Fix $n\geq 2.$ We denote by $L^{2}\left( \mu ^{n}\right) $ the
space of real valued functions on $Z^{n}$ that are
square-integrable with respect to $\mu ^{n}$, and we write
$L_{s}^{2}\left( \mu ^{n}\right) $ to indicate the subspace of
$L^{2}\left( \mu ^{n}\right) $ composed of symmetric functions. %We
%denote by $S_{n}$ and $\widetilde{S}_{n}$, respectively, the
%vector space
%generated by simple functions with the form%
%\begin{equation}
%f(z_{1},\ldots ,z_{n})=\mathbf{1}_{B_{1}}(z_{1})\ldots \mathbf{1}%
%_{B_{n}}(z_{n})\text{,}  \label{sn}
%\end{equation}%
%where $B_{1},\ldots ,B_{n}$ are disjoint subsets of $Z$, and the
%vector space generated by the functions obtained as the
%symmetrization of some element of $S_{n}$. If $f$ is as in
%(\ref{sn}) and $\tilde{f}\in \widetilde{S}_{n}$ is its canonical
%symmetrization,
%we define $I_{n}(\tilde{f})\ $as%
%\begin{equation}
%I_{n}(\tilde{f})=\widehat{N}(B_{1})\cdot \cdot \cdot
%\widehat{N}(B_{n}). \label{in}
%\end{equation}%
For every $f\in L_{s}^{2}(\mu ^{n})$, we denote by $I_{n}(f)$ the
\textsl{multiple Wiener-It\^{o} integral} of order $n$, of $f$ with
respect to $\widehat{N}$. Observe that, for every $m,n\geq 2$, $f\in
L_{s}^{2}(\mu ^{n})$ and $g\in L_{s}^{2}(\mu ^{m})$, one has the
isometric formula (see e.g. \cite{Surg1984}):
\begin{equation}
\mathbb{E}\big[I_{n}(f)I_{m}(g)\big]=n!\langle f,g%
\rangle _{L^{2}(\mu ^{n})}\mathbf{1}_{(n=m)}.  \label{iso}
\end{equation}%
Now fix $n\geq 2$ and $f \in L^{2}\left( \mu ^{n}\right)$ (not
necessarily symmetric) and denote by $\tilde{f}$ the canonical
symmetrization of $f$: for future use we stress that, as a
consequence of Jensen inequality,
\begin{equation}\label{SymInequality}
\|\tilde{f}\|_{L^2{(\mu^n)}}\leq \|f\|_{L^2{(\mu^n)}}.
\end{equation}
The Hilbert space of random variables of the type $I_n(f)$, where
$n\geq 1$ and $f\in L^2_s(\mu^n)$ is called the $n$th \textsl{Wiener
chaos}
associated with $\widehat{N}$. We also use the following conventional notation: $I_{1}\left( f\right) =%
\widehat{N}\left( f\right) $, $f\in L^{2}\left( \mu \right) $;
$I_{n}\left( f\right) =I_{n}(\tilde{f})$, $f\in L^{2}\left( \mu
^{n}\right) $, $n\geq 2$ (this convention extends the definition of
$I_{n}\left( f\right) $ to non-symmetric functions $f$);
$I_{0}\left( c\right) =c$, $c\in \mathbb{R}$. The following
proposition, whose content is known as the \textsl{chaotic
representation property} of $\widehat{N}$, is one of the crucial
results used in this paper. See e.g. \cite{NualartVives} or
\cite{Sur} for a proof.
\begin{prop}[Chaotic decomposition] \label{P : Chaos}
Every random variable $F\in
L^2(\mathcal{F},\mathbb{P})=L^2(\mathbb{P})$ admits a (unique)
chaotic decomposition of the type
\begin{equation}\label{chaosexpansion}
F = \mathbb{E}(F) + \sum_{n\geq 1}^{\infty} I_n(f_n),
\end{equation}
where the series converges in $L^2$ and, for each $n\geq 1$, the
kernel $f_n$ is an element of $L_{s}^{2}(\mu ^{n})$.
\end{prop}
\subsection{Contractions, stars and products}\label{SS : stars}
We now recall a useful version of the \textit{multiplication
formula }for multiple Poisson integrals. To this
end, we define, for $q,p\geq 1$, $f\in L_{s}^{2}\left( \mu ^{p}\right) $, $%
g\in L_{s}^{2}\left( \mu ^{q}\right) $, $r=0,...,q\wedge p$ and $l=1,...,r$,
the (contraction) kernel on $Z^{p+q-r-l}$, which reduces the number of
variables in the product $fg$ from $p+q$ to $p+q-r-l$ as follows: $r$
variables are identified and, among these, $l$ are integrated out. This
contraction kernel is formally defined as follows:
\begin{eqnarray}
&&f\star _{r}^{l}g(\gamma _{1},\ldots ,\gamma _{r-l},t_{1},\ldots
,t_{p-r},s_{1},\ldots ,s_{q-r}) \label{Beirut} \\
&=&\!\!\!\!\int_{Z^{l}}\!\!\!f(z_{1},\ldots ,z_{l},\gamma
_{1},\ldots ,\gamma _{r-l},t_{1},\ldots ,t_{p-r})g(z_{1},\ldots
,z_{l},\gamma _{1},\ldots ,\gamma _{r-l},s_{1},\ldots ,s_{q-r})\mu
^{l}\left( dz_{1}...dz_{l}\right) \text{,} \notag
\end{eqnarray}%
and, for $l=0$,
\begin{equation}
f\star _{r}^{0}g(\gamma _{1},\ldots ,\gamma _{r},t_{1},\ldots
,t_{p-r},s_{1},\ldots ,s_{q-r})=f(\gamma _{1},\ldots ,\gamma
_{r},t_{1},\ldots ,t_{p-r})g(\gamma _{1},\ldots ,\gamma _{r},s_{1},\ldots
,s_{q-r}),  \label{Ocont}
\end{equation}%
so that $f\star _{0}^{0}g(t_{1},\ldots ,t_{p},s_{1},\ldots
,s_{q})=f(t_{1},\ldots ,t_{p})g(s_{1},\ldots ,s_{q})$. For example, if $%
p=q=2 $,
\begin{eqnarray}
&& f\star _{1}^{0}g\left( \gamma ,t,s\right) =f\left( \gamma
,t\right) g\left( \gamma ,s\right) \text{, \ \ }f\star
_{1}^{1}g\left( t,s\right) =\int_{Z}f\left( z,t\right) g\left(
z,s\right) \mu \left( dz\right)
\label{exCont1} \\
&& f\star _{2}^{1}g\left( \gamma \right) =\int_{Z}f\left( z,\gamma
\right) g\left( z,\gamma \right) \mu \left( dz\right) \label{exCont2} \\
&& f\star _{2}^{2}g=\int_{Z}\int_{Z}f\left( z_{1},z_{2}\right)
g\left( z_{1},z_{2}\right) \mu \left( dz_{1}\right) \mu \left(
dz_{2}\right). \label{exCont3}
\end{eqnarray}

\noindent The quite suggestive `star-type' notation is standard,
and has been first used by Kabanov in \cite{Kab} (but see also
Surgailis \cite{Surg1984}). Plainly, for some choice of $f,g,r,l$
the contraction $f\star^l_r g$ may not exist, in the sense that
its definition involves integrals that are not well-defined. On
the positive side, the contractions of the following three types
are well-defined (although possibly infinite) for every $q \geq 2$
and every kernel $f\in L_s^2(\mu^q)$:
\begin{itemize}
\item[(a)] $f\star_r^0 f(z_1,....,z_{2q-r})$, where $r=0,....,q$, as obtained from (\ref{Ocont}), by setting $g=f$;
\item[(b)] $f\star_q^l f(z_1,...,z_{q-l}) = \int_{Z^l} f^2(z_1,...,z_{q-l},\cdot)d\mu^l $, for every
$l=1,...,q$;
\item[(c)] $f\star_r^r f$, for $r=1,....,q-1$.
\end{itemize}

\noindent In particular, a contraction of the type $f\star_q^l f$,
where $l=1,...,q-1$ may equal $+\infty$ at some point
$(z_1,...,z_{q-l})$. The following (elementary) statement ensures
that any kernel of the type $f\star_r^r g$ is square-integrable.

\begin{lemme}\label{L : star-integrability}
Let $p,q\geq 1$, and let $f\in L_s^2(\mu^q)$ and $g\in
L_s^2(\mu^p)$. Fix $r=0,...,q\wedge p$. Then, $f\star_r^r g \in
L^2(\mu^{p+q-2r})$.
\end{lemme}
{\bf Proof.} Just use equation (\ref{Beirut}) in the case $l=r$,
and deduce the conclusion by a standard use of the Cauchy-Schwarz
inequality. \fin

We also record the following identity (which is easily verified by
a standard Fubini argument), valid for every $q\geq 1$, every
$p=1,...,q$:
\begin{equation}\label{useful}
\int_{Z^{2q-p}} (f\star_p^0 f)^2 d\mu^{2q-p} = \int_{Z^p}
(f\star_q^{q-p} f) ^2 d\mu^p,
\end{equation}
for every $f\in L_s^2(\mu^q)$. The forthcoming product formula for
two Poisson multiple integrals is proved e.g. in \cite{Kab} and
\cite{Surg1984}.

\begin{prop}[Product formula]\label{P : ProductPoisson}
Let $f\in L_{s}^{2}\left( \mu ^{p}\right) $ and $g\in
L_{s}^{2}\left( \mu ^{q}\right) $, $p,q\geq 1$, and
suppose moreover that $f\star _{r}^{l}g\in L^{2}(\mu ^{p+q-r-l})$ for every $%
r=0,...,p\wedge q$ and $l=1,...,r$ such that $l\neq r$. Then,
\begin{equation}
I_{p}(f)I_{q}(g)=\sum_{r=0}^{p\wedge q}r!\dbinom{p}{r}\dbinom{q}{r}%
\sum_{l=0}^{r}\binom{r}{l}I_{q+p-r-l}(\widetilde{f\star _{r}^{l}g}),
\label{pproduct}
\end{equation}%
where the tilde `` $\widetilde{}$ '' stands for symmetrization,
that is,
\begin{equation*}
\widetilde{f\star _{r}^{l}g}\left( x_{1},...,x_{q+p-r-l}\right) =\frac{1}{%
\left( q+p-r-l\right) !}\sum_{\sigma }f\star _{r}^{l}g\left( x_{\sigma
\left( 1\right) },...,x_{\sigma \left( q+p-r-l\right) }\right) \text{,}
\end{equation*}%
where $\sigma $ runs over all $\left( q+p-r-l\right) !$ permutations of the
set $\left\{ 1,...,q+p-r-l\right\} $.
\end{prop}

\begin{rem}\label{REM : LEVY}
{\rm  ({\it Multiple integrals and L\'{e}vy processes}). In the multiple Wiener-It\^{o} integrals introduced in this section, the integrators are
compensated Poisson measures of the type $\widehat{N}(dz)$, defined on some abstract Borel space. It is well-known that one can also build similar objects in the framework of L\'{e}vy processes indexed by the real line. Suppose indeed that we are given a cadlag L\'evy process $X=\{X_t,\, t \ge 0\}$ (that is, $X$ has stationary and independent increments, $X$ is
continuous in probability and $X_0=0$), defined on a complete probability space $(\Omega, {\cal F}, P)$,
and with L\'evy triplet $(\gamma,\sigma^2,\nu)$ where $\gamma \in\R,\,
\sigma\ge 0$ and $\nu$ is a L\'evy measure on $\R$ (see e.g.
\cite{Sato}). The process $X$ admits a L\'evy-It{\^o} representation
\begin{equation}
\label{L\'evy-ito} X_t=\gamma t+\sigma W_t+\iint_{(0,t]\times\{\vert
x\vert > 1\}}x\, dN(s,x)+ \lim_{\varepsilon\downarrow
0}\iint_{(0,t]\times\{\varepsilon<\vert x\vert \le 1\}}x\,
 d\widetilde N(s,x),
\end{equation}
 where: (i) $\{W_t,\, t\ge 0\}$ is a standard Brownian motion, (ii)
 $$N(B)=\#\{t:\, (t,\Delta X_t)\in B\},\quad B\in{\cal B}((0,\infty)\times \R_0),$$
 is the jump measure associated with $X$ (where $\R_0=\R-\{0\}$,  $\Delta X_t=X_t-X_{t-}$  and $\#A$ denotes the cardinal of a set $A$), and (iii)
 $$d\widetilde N (t,x)=dN(t,x)-dt\, d\nu(x)$$
 is the compensated jump measure. The convergence in (\ref{L\'evy-ito}) is a.s. uniform (in $t$) on every bounded interval. Following It{\^o} \cite{ito}, the process  $X$ can be extended to an independently scattered random measure $M$ on $(\R_+\times \R, {\cal B}(\R_+\times
 \R))$ as follows.  First,
 consider the   measure $d\mu^*(t,x)=\sigma^2\, dt\,d\delta_0(x)+x^2\,dt\,d\nu(x),$
  where $\delta_0$
  is the Dirac measure at point $0$, and $dt$ is the Lebesgue measure on $\R$. This means that, for
 $E\in {\cal B}(\R_+\times\R)$,
 \begin{equation*}
 \label{mesura}
 \mu^*(E)=\sigma^2\int _{E(0)}dt+\iint _{E'} x^2\, dt\,d\nu(x),
 \end{equation*}
 where $E(0)=\{t\in\R_+:\, (t,0)\in E\}$ and $E'=E-\{(t,0)\in E\}$;
this measure is continuous (see It{\^o} \cite{ito},  p. 256). Now, for
$E\in {\cal B}(\R_+\times\R)$ with $\mu^*(E)<\infty,$
 define
 \begin{equation*}
 \label{mesura-aleat}
M(E)=\sigma\int _{E(0)}dW_t+\lim_n\iint _{\{(t,x)\in E: 1/n<\vert
x\vert<n\}} x\,d\widetilde N(t,x),
  \end{equation*}
 (with convergence in $L^2(\Omega)$), that is, $M$ is a centered independently scattered random
  measure such that
  $$E\big[M(E_1)M(E_2)]=\mu^*(E_1\cap E_2),$$
for $E_1,E_2\in {\cal B}(\R_+\times\R)$ with $\mu^*(E_1)<\infty$ and
$\mu^*(E_2)<\infty$.
Now write $L^2_n=L^2((\R_+\times \R)^n, {\cal B}(\R_+\times \R)^n,
\mu^{*\otimes n}).$
For every $f\in L^2_n$, one can now define a multiple stochastic integral
$I_n(f)$, with integrator $M$, by using the same procedure as in the Poisson case (the first construction of this type is due to It{\^o} -- see \cite{ito}). Since $M$ is defined on a product space, in this scenario we can separate the time and the jump sizes. If  the L\'evy process has no Brownian components, then $M$ has the representation
$M(dz)=M(ds,dx)=x\widetilde{N}(ds,dx)$, that is, $M$ is obtained by integrating a factor $x$
with respect to the underlying compensated Poisson measure of intensity $dt\, \nu(dx)$.
Therefore, in order to define multiple integrals of order $n$ with respect to $M$, one should consider kernels that are square-integrable with respects the measure $(x^2dt\, \nu(dx))^{\otimes n}$. The product formula of multiple integrals in this different context would be analogous to the one given in Proposition \ref{P : ProductPoisson}, but in the definitions of the contraction kernels one has to to take into account the supplementary factor $x$ (see e.g. \cite{lee1} for more details on this point). }
\end{rem}

\subsection{Four Malliavin-type operators} \label{SS : 4Operators}
In this section, we introduce four operators of the
Malliavin-type, that are involved in the estimates of the
subsequent sections. Each of these operators is defined in terms
of the chaotic expansions of the elements in its domain, thus
implicitly exploiting the fact that the chaotic representation
property (\ref{chaosexpansion}) induces an isomorphism between
$L^2(\mathbb{P})$ and the symmetric Fock space canonically
associated with $L^{2}\left( \mu \right)$. The reader is referred
to  \cite{NualartVives} or \cite{SolUtzViv} for more details on
the construction of these operators, as well as for further
relations with chaotic expansions and Fock spaces. One crucial
fact in our analysis is that the derivative operator $D$ (to be
formally introduced below) admits a neat characterization in terms
of a usual difference operator (see the forthcoming Lemma \ref{L :
DifferenceDerivative}).

For later reference, we recall that the space
$L^2(\mathbb{P};L^2(\mu)) \simeq L^2(\Omega\times Z,
\mathcal{F}\otimes \mathcal{Z}, P\otimes \mu)$ is the space of the
measurable random functions $u : \Omega\times Z \rightarrow \R $
such that
$$
\mathbb{E}\left[\int_Z u^2_z \, \mu(dz)\right] <\infty.
$$
In what follows, given $f\in L^2_s(\mu^q)$ ($q\geq 2$) and $z\in Z$,
we write $f(z,\cdot)$ to indicate the function on $Z^{q-1}$ given by
$(z_1,...,z_{q-1})\rightarrow f(z,z_1,...,z_{q-1})$.
\medskip

\noindent i) \underline{The derivative operator $D$}. The
derivative operator, denoted by $D$, transforms random variables
into random functions. Formally, the domain of $D$, written ${\rm
dom} D$, is the set of those random variables $F \in
L^2(\mathbb{P})$ admitting a chaotic decomposition
(\ref{chaosexpansion}) such that
\begin{equation}\label{condDer}
\sum_{n\geq 1}n n! \|f_n\|_{L^2(\mu^n)}^2 < \infty.
\end{equation}
If $F$ verifies (\ref{condDer}) (that is, if $F \in {\rm dom}D$),
then the random function  $z\rightarrow D_zF$ is given by
\begin{equation}\label{TheDerivative}
D_z F = \sum_{n \geq 1}n I_{n-1}(f(z,\cdot)), \,\,\, z\in Z.
\end{equation}
For instance, if $F=I_1(f)$, then $D_z F$ is the non-random
function $z \rightarrow f(z)$. If $F=I_2(f)$, then $D_z F$ is the
random function
\begin{equation}\label{Der2exx}
z\rightarrow 2 I_1(f(z,\cdot)).
\end{equation}
By exploiting the isometric properties of multiple integrals, and
thanks to (\ref{condDer}), one sees immediately that $DF \in
L^2(\mathbb{P};L^2(\mu))$, for every $F \in {\rm dom}D$. Now
recall that, in this paper, the underlying probability space
$(\Omega, \mathcal{F},\mathbb{P})$ is such that $\Omega$ is the
collection of discrete measures given by (\ref{NV1}). Fix $z \in
Z$; given a random variable $F : \Omega \rightarrow \R$, we define
$F_z$ to be the random variable obtained by adding the Dirac mass
$\delta_z$ to the argument of $F$. Formally, for every $\omega \in
\Omega$, we set
\begin{equation}\label{Translation}
F_z(\omega) = F(\omega + \delta_z).
\end{equation}
A version of the following result, whose content is one of the
staples of our analysis, is proved e.g. in \cite[Theorem
6.2]{NualartVives} (see also \cite[Proposition 5.1]{SolUtzViv}). One
should also note that the proof given in \cite{NualartVives} applies
to the framework of a Radon control measure $\mu$: however, the
arguments extend immediately to the case of a Borel control measure
considered in the present paper, and we do not reproduce them here
(see, however, the remarks below). It relates $F_z$ to $D_z F$ and
provides a representation of $D$ as a difference operator.
\begin{lemme}\label{L : DifferenceDerivative}
For every $F \in {\rm dom}D$,
\begin{equation}\label{Der = Diff}
D_z F = F_z - F, \,\,\, {\rm a.e.-}\mu(dz).
\end{equation}
\end{lemme}

\smallskip

\begin{rem} {\rm
\begin{itemize}
\item[1.] In the It\^{o}-type framework detailed in Remark \ref{REM : LEVY},
we could alternatively use the definition of the Malliavin
derivative introduced in \cite{SolUtzViv}, that is,
$$D_{(t,x)} F = \frac{F_{(t,x)} - F}{x}, \,\,\, {\rm
a.e.-} x^2dt\nu(dx).$$ See also \cite{LeSoUtVi}.
\item[2.] When applied to the case $F = I_1(f)$, formula (\ref{Der =
Diff}) is just a consequence of the straightforward relation
$$
F_z -F = \int_Z f(x) (\widehat{N}(dx)+ \delta_z(dx)) - \int_Z f(x)
\widehat{N}(dx) = \int_Z f(x) \delta_z(dx) = f(z).
$$
\item[3.] To have an intuition of the proof of (\ref{Der = Diff}) in the
general case, consider for instance a random variable of the type
$F= \widehat{N}(A)\times \widehat{N}(B)$, where the sets $A,B \in
\mathcal{Z}_\mu$ are disjoint. Then, one has that $F = I_2(f)$,
where $f(x,y) = 2^{-1}\{{\bf 1}_A(x){\bf 1}_B(y) + {\bf 1}_A(y){\bf
1}_B(x)\}$. Using (\ref{Der2exx}), we have that
\begin{equation}\label{A_A}
D_z F = \widehat{N}(B){\bf 1}_A(z) + \widehat{N}(A){\bf 1}_B (z),
\,\, z\in Z,
\end{equation}
and
\begin{equation}F_z -F = \{\widehat{N}+\delta_z\} (A) \times
\{\widehat{N}+\delta_z\}(B) - \widehat{N}(A)\widehat{N}(B), \,\,
z\in Z.\notag
\end{equation}
Now fix $z\in Z$. There are three
possible cases: (i) $z \notin A$ and $z\notin B$, (ii) $z\in A$, and
(iii) $z\in B$. If $z$ is as in (i), then $F_z = F$. If $z$ is as in
(ii) (resp. as in (iii)), then
$$F_z -F = \{\widehat{N}+\delta_z\} (A) \times \widehat{N}(B) -
\widehat{N}(A)\widehat{N}(B)= \widehat{N}(B)$$ (resp.
 $$F_z -F = \widehat{N}(A)\times \{\widehat{N}+\delta_z\}(B) -
\widehat{N}(A)\widehat{N}(B) = \widehat{N}(A)\,\, \,).
 $$
As a consequence, one deduces that $F_z -F = \widehat{N}(B){\bf
1}_A(z) + \widehat{N}(A){\bf 1}_B (z)$, so that relation (\ref{Der
= Diff}) is obtained from (\ref{A_A}).
\item[4.] Observe also that Lemma \ref{L : DifferenceDerivative} yields
that, if $F,G \in {\rm dom}D$ are such that $FG \in {\rm dom}D$,
then
\begin{equation}\label{DerProduct}
D(FG) = FDG + GDF +DGDF,
\end{equation}
(see \cite[Lemma 6.1]{NualartVives} for a detailed proof of this
fact).
\end{itemize}
}
\end{rem}
\medskip

\noindent ii) \underline{The Skorohod integral $\delta$}. Observe
that, due to the chaotic representation property of $\widehat{N}$,
every random function $u \in L^2(\mathbb{P},L^2(\mu)) $ admits a
(unique) representation of the type
\begin{equation}\label{processes}
u_z = \sum_{n\geq 0} I_n(f_n(z,\cdot)),\,\,\,z\in Z,
\end{equation}
where, for every $z$, the kernel $f_n(z,\cdot)$ is an element of
$L^2_s(\mu^n)$. The domain of the Skorohod integral operator,
denoted by ${\rm dom}\delta$, is defined as the collections of
those $u\in L^2(\mathbb{P},L^2(\mu))$ such that the chaotic
expansion (\ref{processes}) verifies the condition
\begin{equation}\label{DomDelta}
\sum_{n\geq 0}(n+1)! \|f_n\|^2_{L^2(\mu^{n+1})}<\infty.
\end{equation}
If $u \in {\rm dom}\delta$, then the random variable $\delta (u)$
is defined as
\begin{equation}\label{Skoro}
\delta(u) = \sum_{n \geq 0} I_{n+1}(\tilde{f}_n),
\end{equation}
where $\tilde{f}_n$ stands for the canonical symmetrization of $f$
(as a function in $n+1$ variables). For instance, if $u(z) = f(z)$
is a deterministic element of $L^2(\mu)$, then $\delta(u) =
I_1(f)$. If $u(z) = I_1(f(z,\cdot))$, with $f\in L^2_s(\mu^2)$,
then $\delta(u)=I_2(f)$ (we stress that we have assumed $f$ to be
symmetric). The following classic result, proved e.g. in
\cite{NualartVives}, provides a characterization of $\delta$ as
the adjoint of the derivative $D$.
\begin{lemme}[Integration by parts formula]\label{L : SkorohodAsAdjoint} For every $G \in {\rm
dom}D$ and every $u \in {\rm dom} \delta$, one has that
\begin{equation}\label{dualDdelta}
\mathbb{E} [G \delta(u)] = \mathbb{E}[\langle DG , u   \rangle _{L^2
(\mu)}],
\end{equation}
where
$$
\langle DG, u  \rangle _{L^2 (\mu)} = \int_Z D_z G \times u(z)
\mu(dz).
$$
\end{lemme}

\medskip

\noindent iii) \underline{The Ornstein-Uhlenbeck generator $L$}. The
domain of the Ornstein-Uhlenbeck generator (see \cite[Chapter
1]{NuaBook}), written ${\rm dom}L$, is given by those $F \in
L^2(\mathbb{P})$ such that their chaotic expansion
(\ref{chaosexpansion}) verifies
$$
\sum_{n\geq 1}n^2 n! \|f_n\|^2_{L^2(\mu^{n})}<\infty.
$$
If $F \in {\rm dom}L $, then the random variable $LF$ is given by
\begin{equation}\label{LF}
LF = -\sum_{n\geq 1} nI_n(f_n).
\end{equation}

Note that $\mathbb{E}(LF)=0$, by definition. The following result
is a direct consequence of the definitions of $D$, $\delta$ and
$L$.
\begin{lemme}\label{L : deltaDL} For every $F \in {\rm dom}L$,
one has that $F \in {\rm dom}D$ and $D F \in{\rm dom}\delta$.
Moreover,
\begin{equation}\label{DdL}
\delta D F = -LF.
\end{equation}
\end{lemme}
{\bf Proof.} The first part of the statement is easily proved by
applying the definitions of ${\rm dom}D$ and ${\rm dom}\delta$
given above. In view of Proposition \ref{P : Chaos}, it is now
enough to prove (\ref{DdL}) for a random variable of the type
$F=I_q(f)$, $q\geq 1$. In this case, one has immediately that $D_z
F = qI_{q-1}(f(z,\cdot))$, so that $\delta D F = qI_q(f) = -LF$.
\fin

\noindent iv) \underline{The inverse of $L$}. The domain of
$L^{-1}$, denoted by $L_0^2(\mathbb{P})$, is the space of
\textsl{centered} random variables in $L^2(\mathbb{P})$. If $F \in
L_0^2(\mathbb{P})$ and $F = \sum_{n\geq 1}I_n(f_n)$ (and thus is centered) then
\begin{equation}\label{L-1F}
L^{-1}F = -\sum_{n\geq 1} \frac1n \, I_n(f_n).
\end{equation}

\subsection{Normal approximation in the Wasserstein distance, via Stein's method}
We shall now give a short account of Stein's method, as applied to
normal approximations in the Wasserstein distance. We denote by
${\rm Lip}(1)$ the class of real-valued Lipschitz functions, from
$\R$ to $\R$, with Lipschitz constant less or equal to one, that is, functions $h$ that are absolutely continuous and satisfy the relation $\|h'\|_{\infty}\leq 1$. Given
two real-valued random variables $U$ and $Y$, the
\textsl{Wasserstein distance} between the laws of $U$ and $Y$,
written $d_W(U,Y)$ is defined as
\begin{equation}\label{Wass}
d_W(U,Y) = \sup_{h\in {\rm Lip}(1)} | \mathbb{E}[h(U)]-
\mathbb{E}[h(Y)] |.
\end{equation}
We recall that the topology induced by $d_W$ on the class of
probability measures over $\R$ is strictly stronger than the
topology of weak convergence (see e.g. \cite{Dudely}). We will
denote by $\mathscr{N}(0,1)$ the law of a centered standard
Gaussian random variable.

Now let $X\sim \mathscr{N}(0,1)$. Consider a real-valued function
$h:\R\rightarrow\R$ such that the expectation $\mathbb{E}[h(X)]$
is well-defined. The \textsl{Stein equation} associated with $h$
and $X$ is classically given by
\begin{equation}\label{SteinGaussEq}
h(x)-\mathbb{E}[h(X)] = f'(x)-xf(x), \quad x\in\R.
\end{equation}
A solution to (\ref{SteinGaussEq}) is a function $f$ depending on
$h$ which is Lebesgue a.e.-differentiable, and such that there
exists a version of $f'$ verifying (\ref{SteinGaussEq}) for every
$x\in\R$. The following result is basically due to Stein
\cite{Stein_orig, Stein_book}. The proof of point (i) (whose content
is usually referred as \textsl{Stein's lemma}) involves a standard
use of the Fubini theorem (see e.g. \cite[Lemma
2.1]{Chen_Shao_sur}). Point (ii) is proved e.g. in \cite[Lemma
4.3]{Chatterjee_AOP}.
\begin{lemme}\label{Stein_Lemma_Gauss}
\begin{enumerate}
\item[\rm (i)] Let $W$ be a random variable. Then, $W\stackrel{\rm
Law}{=}X\sim \mathscr{N}(0,1)$ if, and only if,
\begin{equation}\label{SteinTruc}
\mathbb{E}[f'(W)-Wf(W)]=0,
\end{equation}
for every continuous and piecewise continuously differentiable
function $f$ verifying the relation $\mathbb{E}|f'(X)|$ $<$
$\infty$.
\item[\rm (ii)] If $h$ is absolutely continuous with bounded derivative,
then (\ref{SteinGaussEq}) has a solution $f_h$ which is twice
differentiable and such that $\|f_h'\|_\infty \leq \|h'\|_\infty$
and $\|f_h''\|_\infty \leq 2 \|h'\|_\infty$.
\end{enumerate}
\end{lemme}

Let $\mathscr{F}_{ W}$ denote the class of twice differentiable
functions, whose first derivative is bounded by 1 and whose second
derivative is bounded by 2. If $h\in{\rm Lip}(1)$ and thus
$\|h'\|_{\infty}\leq 1$, then the solution $f_h$, appearing in Lemma
\ref{Stein_Lemma_Gauss}(ii), satisfies $\|f_h'\|_{\infty}\leq 1$ and
$\|f_h''\|_{\infty}\leq 2$, and hence $f_h\in\mathscr{F}_W$.

Now consider a Gaussian random variable $X\sim \mathscr{N}(0,1)$, and let $Y$ be a
generic random variable such that $\mathbb{E}Y^2<\infty$.  By
integrating both sides of (\ref{SteinGaussEq}) with respect to the
law of $Y$ and by using (\ref{Wass}), one sees immediately that
Point (ii) of Lemma \ref{Stein_Lemma_Gauss} implies that
\begin{equation}
d_W(Y,X)\leq \sup_{f\in\mathscr{F}_W}|\mathbb{E}(f'(Y)-Yf(Y))|.
\label{BoundWass}
\end{equation}
Note that the square-integrability of $Y$ implies that the quantity
$\mathbb{E}[Yf(Y)]$ is well defined (recall that $f$ is Lipschitz).
In the subsequent sections, we will show that one can effectively
bound the quantity appearing on the RHS of (\ref{BoundWass}) by
means of the operators introduced in Section \ref{SS : 4Operators}.

\section{A general inequality involving Poisson functionals} \label{S : General Ineq}
\setcounter{equation}{0}
The following estimate, involving the normal approximation of
smooth functionals of $\widehat{N}$, will be crucial for the rest
of the paper. We use the notation introduced in the previous
section.

\begin{thm}\label{T : MainUpperBound}
Let $F\in{\rm dom}D$ be such that $\mathbb{E}(F)=0$. Let $X \sim \mathscr{N}(0,1)$. Then,
\begin{eqnarray}\label{GenUpBound0}
&& d_W(F,X) \leq \mathbb{E}\left[|1-\langle DF, -DL^{-1}F\rangle
_{L^2(\mu)}| \right] + \int_Z \mathbb{E}\left[|D_z F|^2 |D_z
L^{-1}F|\right]\mu(dz) \\
&& \leq  \sqrt{\mathbb{E}\left[(1-\langle DF, -DL^{-1}F\rangle
_{L^2(\mu)})^2\right]} + \int_Z \mathbb{E}\left[|D_z F|^2 |D_z
L^{-1}F|\right]\mu(dz), \label{GenUpBound}
\end{eqnarray}
where we used the standard notation
$$
\langle DF, -DL^{-1}F\rangle _{L^2(\mu)} = - \int_Z [ D_z F \times
D_z L^{-1}F\, ] \,\, \mu (dz).
$$
Moreover, if $F$ has the form $F=I_q(f)$, where $q\geq 1$ and
$f\in L^2_s(\mu^q)$, then
\begin{eqnarray}\label{normInt}
\langle DF, -DL^{-1}F\rangle _{L^2(\mu)} &=
&q^{-1}\|DF\|^2_{L^2(\mu)}
\\ \label{normInt2} \int_Z \mathbb{E}\left[|D_z F|^2 |D_z
L^{-1}F|\right]\mu(dz) &=&q^{-1} \int_Z \mathbb{E}\left[|D_z
F|^3\right]\mu(dz).
\end{eqnarray}
\end{thm}
{\bf Proof.} By virtue of the Stein-type bound (\ref{BoundWass}), it
is sufficient to prove that, for every function $f$ such that
$\|f'\|_\infty \leq 1 $ and $\|f''\|_\infty \leq 2$ (that is, for
every $f \in \mathscr{F}_W$), the quantity $|\mathbb{E}[f'(F) -
Ff(F)]|$ is smaller than the RHS of (\ref{GenUpBound}). To see this,
fix $f \in \mathscr{F}_W$ and observe that, by using
(\ref{Translation}), for every $\omega$ and every $z$ one has that
$D_z f(F)(\omega)=f(F)_z(\omega) -f(F)(\omega) =f(F_z)(\omega)
-f(F)(\omega)$. Now use twice Lemma \ref{L : DifferenceDerivative},
combined with a standard Taylor expansion, and write
\begin{eqnarray} \label{Winter}
D_z f(F) &=& f(F)_z - f(F) =f(F_z)-f(F)  \\
 &=& f'(F)(F_z -F) +R(F_z - F) =f'(F)(D_z F) +R(D_z F), \notag
\end{eqnarray}
where (due to the fact that $\|f''\|_\infty \leq 2$) the mapping $y
\rightarrow R(y)$ is such that $|R(y)|\leq y^2$ for every $y\in\R$.
We can therefore apply (in order) Lemma \ref{L : deltaDL} and Lemma
\ref{L : SkorohodAsAdjoint} (in the case $u = DL^{-1}F$ and $G =
f(F)$) and infer that
$$
\mathbb{E}[Ff(F)] = \mathbb{E}[LL^{-1}Ff(F)]=\mathbb{E}[-\delta
(DL^{-1}F)f(F)] = \mathbb{E}[\langle Df(F),
-DL^{-1}F\rangle_{L^2(\mu)}].
$$
According to (\ref{Winter}), one has that
$$
\mathbb{E}[\langle Df(F), -DL^{-1}F\rangle_{L^2(\mu)}] =
\mathbb{E}[f'(F)\langle DF, -DL^{-1}F\rangle_{L^2(\mu)}] +
\mathbb{E}[\langle R(DF), -DL^{-1}F\rangle_{L^2(\mu)}].
$$
It follows that
$$
|\mathbb{E}[f'(F) - Ff(F)]| \leq |\mathbb{E}[f'(F)(1-\langle DF,
-DL^{-1}F\rangle_{L^2(\mu)})]| + |\mathbb{E}[\langle R(DF),
-DL^{-1}F\rangle_{L^2(\mu)}]|.
$$
By the fact that $\|f'\|_\infty \leq 1$ and by Cauchy-Schwarz,
\begin{eqnarray*}
|\mathbb{E}[f'(F)(1-\langle DF, -DL^{-1}F\rangle_{L^2(\mu)})]| &
\leq & \mathbb{E}[|1-\langle DF, -DL^{-1}F\rangle_{L^2(\mu)}|] \\
& \leq & \sqrt{\mathbb{E}\left[(1-\langle DF, -DL^{-1}F\rangle
_{L^2(\mu)})^2\right]}.
\end{eqnarray*}
On the other hand, one sees immediately that
\begin{eqnarray*}
|\mathbb{E}[\langle R(DF), -DL^{-1}F\rangle_{L^2(\mu)}]| &\leq &
\int_Z \mathbb{E}\left[|R(D_z F)D_z L^{-1}F|\right]\mu(dz)\\
&\leq &\int_Z \mathbb{E}\left[|D_z F|^2 |D_z
L^{-1}F|\right]\mu(dz),
\end{eqnarray*}
Relations (\ref{normInt}) and (\ref{normInt2}) are immediate
consequences of the definition of $L^{-1}$ given in (\ref{L-1F}).
The proof is complete.
 \fin

\begin{rem}
{\rm Note that, in general, the bounds in formulae (\ref{GenUpBound0})--(\ref{GenUpBound}) can be infinite. In the forthcoming Sections \ref{S : MultipleWeak}--\ref{S : EX}, we will exhibit several examples of random variables in ${\rm dom}D$ such that the bounds in the statement of Theorem \ref{T : MainUpperBound} are finite.}
\end{rem}

\begin{rem}
{\rm Let $G$ be an isonormal Gaussian process (see \cite{NuaBook}) over some separable Hilbert space $\H$, and suppose that $F\in L^2(\sigma(G))$ is centered and differentiable in the sense of Malliavin. Then, the Malliavin derivative of $F$, noted $DF$, is a $\H$-valued random element, and in \cite{NouPe07} it is proved that one has the following upper bound on the Wasserstein distance between the law of $F$ and the law of $X\sim \mathscr{N}(0,1)$:
$$
d_W (F, X) \leq \mathbb{E}|1- \langle DF, -DL^{-1}F\rangle_\H |,
$$
where $L^{-1}$ is the inverse of the Ornstein-Uhlenbeck generator associated with $G$.}
\end{rem}

\medskip

Now fix $h \in L^2(\mu)$. It is clear that the random variable
$I_1(h)=\widehat{N}(h)$ is an element of ${\rm dom} D$. One has that
$D_z\widehat{N}(h) = h(z)$; moreover
$-L^{-1}\widehat{N}(h)=\widehat{N}(h)$. Thanks to these relations,
one deduces immediately from Theorem \ref{T : MainUpperBound} the
following refinement of Part A of Theorem 3 in \cite{PecTaq07}.

\begin{cor}\label{P : 1stChaosBound}
Let $h \in L^2(\mu)$ and let $X\sim\mathscr{N}(0,1)$. Then, the
following bound holds:
\begin{equation}\label{1stChUB}
d_W(\widehat{N}(h), X)\leq \left |1 - \|h\|^2_{L^2(\mu)} \right |
+ \int_Z |h(z)|^3 \mu(dz).
\end{equation}
As a consequence, if $\mu(Z) = \infty$ and if a sequence $\{h_k\}
\subset L^2(\mu) \cap L^3(\mu)$ verifies, as $k\rightarrow
\infty$,
\begin{equation}\label{1CONDclt}
\|h_k\|_{L^2(\mu)}\rightarrow 1 \text{\ \ and \ \ }
\|h_k\|_{L^3(\mu)}\rightarrow 0,
\end{equation}
one has the CLT
\begin{equation}\label{1clt}
\widehat{N}(h_k) \stackrel{\rm law}{\longrightarrow} X,
\end{equation}
and the inequality (\ref{1stChUB}) provides an explicit upper
bound in the Wasserstein distance.
\end{cor}

In the following section, we will
use Theorem \ref{T : MainUpperBound} in order to prove general
bounds for the normal approximation of multiple integrals of
arbitrary order. As a preparation, we now present two simple
applications of Corollary \ref{P : 1stChaosBound}.

\begin{example}
{\rm Consider a centered Poisson measure $\widehat{N}$ on $Z =
\R^+$, with control measure $\mu$ equal to the Lebesgue measure.
Then, the random variable $k^{-1/2}\widehat{N}([0,k])
=\widehat{N}(h_k) $, where $h_k =k^{-1/2}{\bf 1}_{[0,k]}$, is an
element of the first Wiener chaos associated with $\widehat{N}$.
Plainly, since the random variables $\widehat{N}((i-1,i])$
($i=1,...,k$) are i.i.d. centered Poisson with unitary variance, a
standard application of the Central Limit Theorem yields that, as
$k\rightarrow \infty$, $\widehat{N}(h_k) \stackrel{\rm
law}{\rightarrow} X \sim \mathscr{N}(0,1)$. Moreover,
$\widehat{N}(h_k) \in {\rm dom} D$ for every $k$, and
$D\widehat{N}(h_k) = h_k$. Since
$$
\|h_k\|^2_{L^2(\mu)} = 1 \text{\ \ and \ \ } \int_Z |h_k|^3\mu(dz)
= \frac{1}{k^{1/2}},
$$
one deduces from Corollary \ref{P : 1stChaosBound} that
$$
d_W(\widehat{N}(h_k), X )\leq \frac{1}{k^{1/2}},
$$
which is consistent with the usual Berry-Esséen estimates.} \fin
\end{example}
\begin{example} \label{EX : OUlinear} {\rm
Fix $\lambda >0$. We consider the\textsl{\ Ornstein-Uhlenbeck
L\'{e}vy process} given by
\begin{equation}
Y_{t}^{\lambda }=\sqrt{2\lambda }\int_{-\infty
}^{t}\int_{\mathbb{R}}u\exp \left( -\lambda \left( t-x\right)
\right) \widehat{N}\left( du,dx\right) \text{, \ \ }t\geq
0\text{,}  \label{OUl}
\end{equation}%
where $\widehat{N}$ is a centered Poisson measure over
$\mathbb{R\times R}$, with control measure given by $\mu (du, dx)
= \nu \left( du\right) dx$, where $\nu \left(
\cdot \right) $ is positive, non-atomic and normalized in such a way that $\int_{\mathbb{%
R}}u^{2}d\nu$ $=1$. We assume also that $\int_{\mathbb{R}%
}\left\vert u\right\vert ^{3}d\nu <\infty .$ Note that $%
Y_{t}^{\lambda }$ is a stationary \textsl{moving average L\'{e}vy process}.
According to \cite[Theorem 4]{PecTaq07}, on has that, as
$T\rightarrow \infty $,
$$
A^{\lambda}_T = \frac{1}{\sqrt{2T/\lambda}}\int_{0}^{T}Y_{t}^{\lambda }dt\overset{\rm law}{\rightarrow }%
X\sim \mathscr{N}\left( 0,1\right).
$$
We shall show that Corollary \ref{P : 1stChaosBound} implies the
existence of a finite constant $q_\lambda >0$, depending uniquely on
$\lambda$, such that, for every $T>0$
\begin{equation}\label{estate}
d_W(A_T^{\lambda}, X) \leq \frac{q_\lambda}{T^{1/2}}.
\end{equation}
To see this, first use a Fubini theorem in order to represent
$A_T^{\lambda}$ as an integral with respect to $\widehat{N}$, that
is, write
\begin{equation*}
A_T^{\lambda}=\int_{-\infty }^{T}\int_{%
\mathbb{R}}u\left[ \left( \frac{2\lambda }{2T/\lambda}\right) ^{\frac{1}{2}%
}\int_{x\vee 0}^{T}\exp \left( -\lambda \left( t-x\right) \right)
dt\right] \widehat{N}\left( du,dx\right):=
\widehat{N}(h^{\lambda}_T).
\end{equation*}
Clearly, $A_T^{\lambda} \in {\rm dom} D$ and $DA_T^{\lambda} =
h_T^{\lambda}$. By using the computations contained in \cite[Proof
of Theorem 4]{PecTaq07}, one sees immediately that there exists a
constant $\beta_\lambda$, depending uniquely on $\lambda$ and such
that
$$
\int_{\R\times\R} |h_T^{\lambda}(u,x)|^3\nu(du)dx \leq
\frac{\beta_\lambda}{T^{1/2}}.
$$
Since one has also that $|\|h_T^{\lambda}\|^2_{L^2(\mu)}-1 |=
O(1/T)$, the estimate (\ref{estate}) is immediately deduced from
Corollary \ref{P : 1stChaosBound}. \fin }
\end{example}

\section{Bounds on the Wasserstein distance for multiple integrals of arbitrary order}\label{S : MultipleWeak}\setcounter{equation}{0} In this section we establish
general upper bounds for multiple Wiener-Itô integrals of arbitrary order
$q\geq 2$, with respect to the compensated Poisson measure
$\widehat{N}$. Our techniques hinge on the forthcoming Theorem \ref{T
: MWII-Bounds}, which uses the product formula (\ref{pproduct})
and the inequality (\ref{GenUpBound}).

\subsection{The operators $G_p^q$ and $\widehat{G}_p^q$}\label{SS : FurNot}
Fix $q\geq 2$ and let $f\in L^2_s(\mu^q)$. The operator $G_p^q$
transforms the function $f$, of $q$ variables, into a function
$G_p^qf$ of $p$ variables, where $p$ can be as large as $2q$. When
$p=0$, we set
$$G_0^q f = q!\|f\|^2_{L^2(\mu^q)}$$ and, for every $p=1,...,2q$, we
define the function $(z_1,...,z_p)\rightarrow G_p^q f(z_1,...,z_p)$,
from $Z^p$ into $\R$, as follows:
\begin{equation}\label{Gpi}
G_p^q f(z_1,...,z_p) = \sum_{r=0}^q \sum_{l=0}^r {\bf 1}_{\{2q-r-l
=p\}} r! \binom{q}{r}^2 \binom{r}{l} \widetilde{f\star_r^l
f}(z_1,...,z_p),
\end{equation}
where the `star' contractions have been defined in formulae
(\ref{Beirut}) and (\ref{Ocont}), and the tilde ``$\, \widetilde{
}\, $
 '' indicates a symmetrization. The notation (\ref{Gpi}) is mainly
introduced in order to give a more compact representation of the
RHS of (\ref{pproduct}) when $g=f$. Indeed, suppose that $f\in
L^2_s(\mu^q)$ ($q\geq 2$), and that
$f\star _{r}^{l}f\in L^{2}(\mu ^{2q-r-l})$ for every $%
r=0,...,q$ and $l=1,...,r$ such that $l\neq r$; then, by using
(\ref{pproduct}) and (\ref{Gpi}), one deduces that
\begin{equation}
I_{q}(f)^2=\sum_{p=0}^{2q} I_p (G_p^q f), \label{CompactProduct}
\end{equation}
where $I_0 (G_0^q f)= G_0^q f=q!\|f\|^2_{L^2(\mu^q)}$.  Note that
the advantage of (\ref{CompactProduct}) (over (\ref{pproduct})) is
that the square $I_{q}(f)^2$ is now represented as an {\sl
orthogonal} sum of multiple integrals. As before, given $f\in
L^2_s(\mu^q)$ ($q\geq 2$) and $z\in Z$, we write $f(z,\cdot)$ to
indicate the function on $Z^{q-1}$ given by
$(z_1,...,z_{q-1})\rightarrow f(z,z_1,...,z_{q-1})$. To simplify
the presentation of the forthcoming results, we now introduce a
further assumption.

\medskip

\noindent\underline{\bf Assumption A}. For the rest of the paper,
whenever considering a kernel $f\in L_s^2(\mu^q)$, we will
implicitly assume that every contraction of the type
$$(z_1,...,z_{2q-r-l})\longrightarrow | f | \star_r^l | f | (z_1,...,z_{2q-r-l})$$ is
\textsl{well-defined and finite} for every $r=1,...,q$, every
$l=1,...,r$ and every $(z_1,...,z_{2q-r-l})\in Z^{2q-r-l}$.

\medskip

Assumption A ensures, in particular, that the following relations
are true for every $r=0,....,q-1$, every $l=0,...,r$ and every
$(z_1,...,z_{2(q-1)-r-l})\in Z^{2(q-1)-r-l}$:

$$
\int_Z [f(z, \cdot) \star_r^l f(z, \cdot)]\mu(dz) \!=\! \{f
\star_{r+1}^{l+1} f\}, \!\!\text{ \ \ \ and \ \ \ } \!\! \int_Z
[\widetilde{f(z, \cdot) \star_r^l f(z, \cdot)}]\mu(dz) \!=\!
\left\{\widetilde{f \star_{r+1}^{l+1} f}\right \}
$$
(note that the symmetrization on the LHS of the second equality
\textsl{does not} involve the variable $z$).

The notation $G_p^{q-1} f(z,\cdot)$ ($p=0,...,2(q-1)$) stands for
the action of the operator $G_p^{q-1}$ (defined according to
(\ref{Gpi})) on $f(z,\cdot)$, that is,
\begin{eqnarray}\label{oto}
&& G_p^{q-1} f(z,\cdot)(z_1,...,z_p) \\
&& = \sum_{r=0}^{q-1} \sum_{l=0}^r {\bf 1}_{\{2(q-1)-r-l =p\}} r!
\binom{q-1}{r}^2 \binom{r}{l} \widetilde{f(z,\cdot)\star_r^l
f(z,\cdot)}(z_1,...,z_p). \notag
\end{eqnarray}
\noindent For instance, if $f\in L^2_s(\mu^3)$, then
\begin{equation*}\notag
G_1^2 f(z,\cdot)(a) = 4 \times f(z,\cdot)\star_2^1 f(z,\cdot)(a) =
4\int_Z f(z,a,u)^2\mu(du).
\end{equation*}
Note also that, for a fixed $z\in Z$, the quantity
$G_0^{q-1}f(z,\cdot)$ is given by the following constant
\begin{equation}\label{G0}
G_0^{q-1}f(z,\cdot) =(q-1)! \int_{Z^{q-1}} f^2(z,\cdot)d\mu^{q-1}.
\end{equation}
Finally, for every $q\geq 2$ and every $f\in L^2_s(\mu^q)$, we set
$$\widehat{G}_0^q f =\int_Z G_0^{q-1}f(z,\cdot)\mu(dz) =
q^{-1}G_0^q f = (q-1)!\|f\|^2_{L^2(\mu^q)},$$ and, for
$p=1,...,2(q-1)$, we define the function $(z_1,...,z_p)\rightarrow
\widehat{G}_p^q f(z_1,...,z_p)$, from $Z^p$ into $\R$, as follows:
\begin{equation}\label{GpiHat}
\widehat{G}_p^q f(\cdot) =\int_Z G_p^{q-1}f(z,\cdot)\mu(dz),
\end{equation}
or, more explicitly,
\begin{eqnarray}
&& \widehat{G}_p^q f(z_1,...,z_p)\notag \\
\label{GpiHat2} && = \int_Z \sum_{r=0}^{q-1} \sum_{l=0}^r {\bf 1}
_{\{2(q-1)-r-l=p\}}
r!\binom{q-1}{r}^2\binom{r}{l} \widetilde{f(z,\cdot) \star_r^l f(z,\cdot)}(z_1,...,z_p)\mu(dz)\\
&&=\sum_{t=1}^{q} \sum_{s=1}^{t\wedge (q-1) } {\bf 1}
_{\{2q-t-s=p\}}
(t-1)!\binom{q-1}{t-1}^2\binom{t-1}{s-1}\widetilde{f\star_t^s
f}(z_1,...,z_p), \label{Gpihat3}
\end{eqnarray}
where in (\ref{Gpihat3}) we have used the change of variables
$t=r+1$ and $s=l+1$, as well as the fact that $p\geq 1$. We stress
that the symmetrization in (\ref{GpiHat2}) does not involve the
variable $z$.
\subsection{Bounds on the Wasserstein distance}\label{SS : UpperBchaos}
When $q\geq 2$ and $\mu(Z) = \infty$, we shall focus on kernels
$f\in L^2_s(\mu^q)$ verifying the following technical condition:
for every $p=1,...,2(q-1)$,
\begin{equation}\label{technical}
\int_Z \left[ \sqrt{\int_{Z^p} \{G_p^{q-1}f(z,\cdot)\}^2 \,\,
d\mu^p}\,\,\,\, \right] \mu(dz) < \infty.
\end{equation}
\noindent As will become clear from the subsequent discussion, the
requirement (\ref{technical}) is used to justify a Fubini argument.

\begin{rem}\label{R : Tech} {\rm
\begin{itemize}
\item[1.] When $q=2$, one deduces from (\ref{oto}) that $G_1^1f(z,\cdot) (x) = f(z,x)^2$ and
$G_2^1 f(z,\cdot)(x,y) = f(z,x)f(z,y)$. It follows that, in this
case, condition (\ref{technical}) is verified if, and only if, the
following relation holds:
\begin{equation}\label{technical_q=2}
\int_Z  \sqrt{\int_{Z} f(z,a)^4 \,\, \mu(da)}\,\,\,\,\mu(dz) <
\infty. \end{equation}

Indeed, since $f$ is square-integrable, the additional relation
\begin{equation}\label{fake}  \int_Z  \sqrt{\int_{Z^2} f(z,a)^2 f(z,b)^2
\mu(da)\mu(db)}\,\,\,\,\mu(dz) < \infty
\end{equation}
is always satisfied, since
$$\int_Z \sqrt{\int_{Z^2} f(z,a)^2 f(z,b)^2 \mu(da)\mu(db)}\mu(dz)
= \|f\|^2_{L^2(\mu^2)}.$$
\item[2.] From relation (\ref{oto}), one deduces immediately that (\ref{technical}) is implied by the following
(stronger) condition: for every $p=1,...,2(q-1)$ and every $(r,l)$
such that $2(q-1)-r-l =p$
\begin{equation}
\int_Z \left[ \sqrt{\int_{Z^p}\{ f(z,\cdot)\star_r^l
f(z,\cdot)\}^2 \,\, d\mu^p}\,\,\,\, \right] \mu(dz) < \infty.
\end{equation}
\item[3.] When $\mu(Z)< \infty$ and
$$\int_Z \left[ \int_{Z^p} \{G_p^{q-1}f(z,\cdot)\}^2 \,\, d\mu^p
\right] \mu(dz) <\infty,$$ condition (\ref{technical}) is
automatically satisfied (to see this, just apply the
Cauchy-Schwarz inequality).
\item[4.] Arguing as in the previous point, a sufficient condition for
(\ref{technical}) to be satisfied, is that the support of the
symmetric function $f$ is contained in a set of the type $A\times
\cdot \cdot\cdot \times A$ where $A$ is such that $\mu(A) <\infty$ .
\end{itemize}
}
\end{rem}

\begin{thm}[Wasserstein bounds on a fixed chaos]\label{T : MWII-Bounds}
Fix $q\geq 2$ and let $X \sim \mathscr{N}(0,1)$. Let $f\in
L^2_s(\mu^q)$ be such that:
\begin{itemize}
\item[\rm ({\bf i})] whenever $\mu(Z) = \infty$, condition (\ref{technical})
is satisfied for every $p=1,...,2(q-1)$;
\item[\rm ({\bf ii})] for $d\mu$-almost every $z\in Z$, every $r=1,...,q-1$ and
every $l=0,...,r-1$, the kernel $f(z,\cdot)\star_{ r}^{ l}
f(z,\cdot)$ is an element of $L^2_s(\mu^{2(q-1)-r-l})$.
\end{itemize}
Denote by $I_q(f)$ the multiple Wiener-Itô integral, of order $q$,
of $f$ with respect to $\widehat{N}$. Then, the following bound
holds:
\begin{eqnarray} && d_W(I_q(f), X) \leq  \notag \\ \label{bound_on_int}
&& \sqrt{(1-q!\|f\|^2_{L^2(\mu^q)})^2 + q^2
\sum_{p=1}^{2(q-1)} p! \int_{Z^p} \left\{\widehat{G}^q_p f\right\}^2 d\mu^p } \\
&&+ \,\,q^2\sqrt{(q-1)!\|f\|^2_{L^2(\mu^q)}}\times
\sqrt{\sum_{p=0}^{2(q-1)}p!\int_Z
\left\{\int_{Z^p}G_p^{q-1}f(z,\cdot)^2 d\mu^p \right\} \mu(dz) },
\label{bound_on_int_cont}
\end{eqnarray}
where the notations $\widehat{G}^q_p f$ and $G_p^{q-1}f(z,\cdot)$
are defined, respectively, in (\ref{GpiHat})--(\ref{Gpihat3}) and
(\ref{oto}). Moreover, the bound appearing on the  RHS of
(\ref{bound_on_int})--(\ref{bound_on_int_cont}) can be assessed by
means of the following estimate:
\begin{eqnarray}\label{starestimate1a}
&& \sqrt{(1-q!\|f\|^2_{L^2(\mu^q)})^2 + q^2 \sum_{p=1}^{2(q-1)} p!
\int_{Z^p} \left\{\widehat{G}^q_p f\right\}^2 d\mu^p } \leq
|1-q!\|f\|^2_{L^2(\mu^q)}| \\
&&  + \,\,q  \sum_{t=1}^q \sum_{s=1}^{t\wedge(q-1)} {\bf 1}_{\{2
\leq t+s \leq 2q-1 \}}  (2q-t-s)!^{1/2} (t-1)!^{1/2}
\times \label{starestimate1b} \\
&& \text{\ \ \ \ \ \ \ \ \ \ \ \ \ \ \ \ \ \ \ \ \ \ \ \ \ \ \ \ \
\ \ \ \ \ \ \ \ \ \ \ \ \ \ \ \ \ \ \ \ } \times
\binom{q-1}{t-1}\binom{t-1}{s-1}^{1/2} \|f\star_t^s f \|_{L^2
(\mu^{2q-t-s})}, \notag
\end{eqnarray}
Also, if one has that
\begin{eqnarray}
&& f\star_q^{q-b} f \in L^2(\mu^{b}), \,\,\forall b=1,...,q-1,
\label{Declan}
\end{eqnarray}
then,
\begin{eqnarray}\label{starestimate2a}
&& \sqrt{\sum_{p=0}^{2(q-1)}p! \int_Z
\left\{\int_{Z^p}G_p^{q-1}f(z,\cdot)^2 d\mu^p \right\}^2 \mu(dz) }
\leq \\ \label{starestimate2b}&& \sum_{b=1}^q \sum_{a=0}^{b-1}
{\bf 1}_{\{1\leq a+b \leq 2q-1 \}} (a+b)!^{1/2} (q-a-1)!^{1/2} \times \\
\notag && \text{\ \ \ \ \ \ \ \ \ \ \ \ \ \ \ \ \ \ \ \ \ \ \ \ \
\ \ \ \ \ \ \ } \times \binom{q-1}{q-1-a}\binom{q-1-a}{q-b}^{1/2}
\|f\star_b^a f \|_{L^2(\mu^{2q-a-b})}.
\end{eqnarray}
\end{thm}

\medskip

\begin{rem}
{\rm
\begin{itemize}
\item[1.] There are contraction norms in (\ref{starestimate2b}) that do not appear in the previous formula
(\ref{starestimate1b}), and vice versa. For example, in
(\ref{starestimate2b}) one has $\|f\star_b^0 f
\|_{L^2(\mu^{2q-b})}$, where $b=1,...,q$ (this corresponds to the
case $b\in \{1,...,q\}$ and $a=0$). By using formula (\ref{useful}),
these norms can be computed as follows:
\begin{eqnarray} \label{revolved2}
 \|f\star_b^0 f \|_{L^2(\mu^{2q-b})} &=& \|f\star_q^{q-b} f
\|_{L^2(\mu^{b})}, \,\,\,\, b=1,...,q-1 ;
\\ \label{revolved3}
\|f\star_q^0 f \|_{L^2(\mu^{q})} &=&\sqrt{ \int_{Z^q} f^4 d\mu^q
}.
\end{eqnarray}
We stress that $f^2 = f\star_q^0 f$, and therefore $f\star_q^0 f \in
L^2(\mu^q)$ if, and only if, $f \in L^4(\mu^q)$.
\item[2.] One should compare Theorem \ref{T : MWII-Bounds} with Proposition
3.2 in \cite{NouPe07}, which provides upper bounds for the normal
approximation of multiple integrals with respect to an {\sl
isonormal Gaussian process}. The bounds in \cite{NouPe07} are also
expressed in terms of contractions of the underlying kernel.
\end{itemize}
}
\end{rem}

\begin{example}(\textsl{Double integrals}). \label{E :
MWII-2-p1} {\rm We consider a double integral of the type $I_2(f)$,
where $f\in L_s^2(\mu^2)$ satisfies (\ref{technical_q=2}) (according
to Remark \ref{R : Tech}(1), this implies that (\ref{technical}) is
satisfied). We suppose that the following three conditions are
satisfied: (a) $\mathbb{E}I_2(f)^2 = 2\|f\|_{L^2(\mu^2)}^2 = 1$, (b)
$f\star_2^1 f\in L^2 (\mu^{1})$ and (c) $f\in L^4(\mu^2)$. Since
$q=2$ here, one has $f(z,\cdot)\star_1 ^0 f(z,\cdot)(a)= f(z,a)^2$
which is square-integrable, and hence Assumption (ii) and condition
(\ref{Declan}) in Theorem \ref{T : MWII-Bounds} are verified. Using
relations (\ref{revolved2})--(\ref{revolved3}), we can deduce the
following bound on the Wasserstein distance between the law of
$I_2(f)$ and the law of $X\sim\mathscr{N}(0,1)$:
\begin{eqnarray}
&& d_W (I_2(f),X) \label{doubleEX} \\
&& \leq   \sqrt{8} \|f\star_1^1 f\|_{L^2(\mu^2)} +
\left\{2+\sqrt{8}(1+\sqrt{6})\right\}\|f\star_2^1 f\|_{L^2(\mu^2)}
+ 4 \sqrt{\int_{Z^2}f^4 d\mu^2}. \notag
\end{eqnarray}
To obtain (\ref{doubleEX}), observe first that, since assumption (a)
above is in order, then relations
(\ref{bound_on_int})--(\ref{bound_on_int_cont}) in the statement of
Theorem \ref{T : MWII-Bounds} yield that
$$ d_W (I_2(f),X) \leq 0 + (\ref{starestimate1b}) + 2^{3/2} \times
(\ref{starestimate2b}).$$ To conclude observe that
$$(\ref{starestimate1b}) =2 \{ 2^{1/2} \|f\star _1^1 f
\|_{L^2(\mu^2)} + \|f\star_2^1 f\|_{L^2(\mu)} \}$$ and
$$(\ref{starestimate2b}) =  \|f\star _2^1 f \|_{L^2(\mu)} \{3!^{1/2} +1
\} + 2^{1/2} \|f\|^2_{L^4(\mu^2)},$$ since $\|f\star _2^1 f
\|_{L^2(\mu)} = \|f\star _1^0 f \|_{L^2(\mu^3)}$. A general
statement, involving random variables of the type $F=I_1(g)+I_2(h)$
is given in Theorem \ref{T : Bound1C2C}. } \fin
\end{example}

\begin{example}(\textsl{Triple integrals}). \label{E :
MWII-3-p1} {\rm
 We consider a random variable of the type $I_3(f)$,
where $f\in L_s^2(\mu^3)$ verifies (\ref{technical}) (for instance,
according to Remark \ref{R : Tech}(4), we may assume that $f$ has
support contained in some rectangle of finite $\mu^3$-measure). We
shall also suppose that the following three conditions are
satisfied: (a) $\mathbb{E}I_3(f)^2 = 3!\|f\|_{L^2(\mu^3)}^2 = 1$,
(b) For every $r=1,...,3$, and every $l=1,...,r \! \wedge \! 2$, one
has that $f \star_r^l f \! \in \! L^2(\mu^{6-r-l})$, and (c) $f\in
L^4 (\mu^3)$. One can check that all the assumptions in the
statement of Theorem \ref{T : MWII-Bounds} (in the case $q=3$) are
satisfied. In view of (\ref{revolved2})--(\ref{revolved3}), we
therefore deduce (exactly as in the previous example and after some
tedious bookkeeping!) the following bound on the Wasserstein
distance between the law of $I_3(f)$ and the law of
$X\sim\mathscr{N}(0,1)$:
\begin{eqnarray}\label{ZZ}
d_W (I_3(f),X) & \leq &  3(4!^{1/2})\|f\star_1^1 f\|_{L^2(\mu^4)}
+ (6+2\sqrt{18})(3!^{1/2})\|f\star_2^1 f\|_{L^2(\mu^3)} \\
&& +6\sqrt{2}\|f\star_2^2 f\|_{L^2(\mu^2)} +
\{6+\sqrt{18}(4+2(4!^{1/2}))\}\|f\star_3^1 f\|_{L^2(\mu^2)} \notag \\
&& +\{\sqrt{18}(2^{1/2}+5!^{1/2})+6\}\|f\star_3^2 f\|_{L^2(\mu^1)}
\notag \\
&& +\{3^{3/2}(18^{1/2})\} \sqrt{\int_{Z^3}f^4 d\mu^3}. \notag
\end{eqnarray}
}\fin
\end{example}

\noindent {\bf Proof of Theorem \ref{T : MWII-Bounds}}. First
observe that, according to Theorem \ref{T : MainUpperBound}, we
have that
$$
d_W(I_q(f),X) \leq \sqrt{\mathbb{E}\left[\left(1-\frac{1}{q}\|
DI_q(f) \|^2_{L^2(\mu)}\right)^2\right]} + \frac{1}{q}\int_Z
\mathbb{E}|D_z I_q(f)|^3 \mu(dz).
$$
The rest of the proof is divided in four steps:
\begin{itemize}
\item[\bf ({S1})] Proof of the fact that $\sqrt{\mathbb{E}\left[\left(1-\frac{1}{q}\| DI_q(f)
\|^2_{L^2(\mu)}\right)^2\right]}$ is less or equal to the quantity
appearing at the line (\ref{bound_on_int}).
\item[\bf ({S2})] Proof of the fact that $\frac{1}{q}\int_Z
\mathbb{E}|D_z I_q(f)|^3 \mu(dz)$ is less or equal to the quantity
appearing at the line (\ref{bound_on_int_cont}).
\item[\bf ({S3})] Proof of the estimate displayed in formulae
(\ref{starestimate1a})--(\ref{starestimate1b}).
\item[\bf ({S4})] Proof of the estimate in formulae
(\ref{starestimate2a})--(\ref{starestimate2b}), under the
assumption (\ref{Declan}).
\end{itemize}
\noindent \underline{\bf Step (S1)}. Start by defining $A_f$ as
the collection of those $z\in Z$ such that Assumption (ii) in the
statement is violated, that is,
$$
A_f = \left \{z\in Z : \exists r=1,...,q-1, \,\, l=0,...,r-1 :
f(z,\cdot)\star_{ r}^{ l} f(z,\cdot) \notin L^2_s(\mu^{2(q-1)-r-l})
\right \}.
$$
By assumption, one has that $\mu(A_f) = 0$. By writing
$$
A_f^q = \underbrace{A_f \times \cdot\cdot\cdot \times A_f }_{q \,\,
{\rm times}} \in \mathcal{Z}^q,
$$
one also deduces that $\mu^q(A_f^q) = 0$. As a consequence,
$I_q(f) = I_q(f {\bf 1}_{A_f^q})$, a.s.-{$\mathbb{P}$}. It follows
that, by replacing $f$ with $ f {\bf 1}_{A_f^q}$ (and without loss
of generality), we can assume for the rest of the proof that the
following stronger assumption is verified:

\begin{itemize}

\item[({\bf ii}$'$)] {\sl For \underline{every} $z\in Z$, every $r=1,...,q-1$ and
every $l=0,...,r-1$, one has that $f(z,\cdot)\star_{ r}^{ l}
f(z,\cdot)\in L^2_s(\mu^{2(q-1)-r-l})$.}

\end{itemize}

\noindent Now use (\ref{TheDerivative}) to write $D_zI_q(f) =
qI_{q-1}(f(z,\cdot))$, $z\in Z$. By virtue of the multiplication
formula (\ref{pproduct}), in particular (\ref{CompactProduct}), and by adopting the notation
(\ref{Gpi}), we infer that, for every $z\in Z$
\begin{equation}\label{P0}
\{D_zI_q(f)\}^2 = q q! \int_{Z^{q-1}} f^2(z,\cdot)d\mu^{q-1} +
q^2\sum_{p=1}^{2(q-1)} I_p(G_p^{q-1}f(z,\cdot)).
\end{equation}
Since (\ref{technical}) is in order, one has that, for every
$p=1,...,2(q-1)$,
\begin{equation}\label{P1}
\mathbb{E} \int_Z |I_{p}(G_p^{q-1}f(z,\cdot))|\mu(dz) \leq \int_Z
\sqrt{p! \int_{Z^p}G_p^{q-1}f(z,\cdot)^2 d\mu^p} \mu(dz) <\infty,
\end{equation}
where we have used the Cauchy-Schwarz inequality, combined with
the isometric properties of multiple integrals, as well as the
fact that, for every $z\in Z$,
$$\int_{Z^p}\left[\widetilde{G_p^{q-1}f(z,\cdot)}\right]^2 d\mu^p
\leq \int_{Z^p}G_p^{q-1}f(z,\cdot)^2 d\mu^p.$$ Relations
(\ref{P0})--(\ref{P1}) yield that one can write
\begin{eqnarray}\label{P2}
\frac{1}{q}\| DI_q(f) \|^2_{L^2(\mu)}-1 &=&  \frac{1}{q} \int_Z
\{D_zI_q(f)\}^2 \mu(dz)-1 \\ &=& q!\|f\|^2_{L^2(\mu^q)}-1 +
q\sum^{2(q-1)}_{p=1}\int_Z I_p(G_p^{q-1}f(z,\cdot)) \mu(dz).
\notag
\end{eqnarray}
Since assumption (\ref{technical}) is in order, one has that, for
every $p=1,...,2(q-1)$,
\begin{eqnarray}
&& \mathbb{E}\left[\left(\int_Z I_p(G_p^{q-1}f(z,\cdot))
\mu(dz)\right)^2 \right] \label{Sco}\\
&& \leq \!\! \int_{Z^2}\!\!\!
\mathbb{E}\left|I_p(G_p^{q-1}f(z,\cdot))I_p(G_p^{q-1}f(z',\cdot))\right|\mu(dz)\mu(dz')
\notag \\
&& \leq p\,! \, \left\{\int_Z \sqrt{\int_{Z^p} \{G_p^{q-1}f(z,\cdot)\}^2
d\mu^p}\,\, \mu(dz) \right\}^2 <\infty, \notag
\end{eqnarray}
and one can easily verify that, for $1\leq p\neq l \leq 2(q-1)$,
the random variables $$\int_Z I_p(G_p^{q-1}f(z,\cdot))
\mu(dz)\text{\ \ and \ \ } \int_Z I_l(G_l^{q-1}f(z,\cdot))
\mu(dz)$$ are orthogonal in $L^2(\mathbb{P})$. It follows that
$$
\mathbb{E}\left[\left(1-\frac{1}{q}\| DI_q(f)
\|^2_{L^2(\mu)}\right)^2\right] = (q!\|f\|^2_{L^2(\mu^q)}-1)^2 +
q^2\sum^{2(q-1)}_{p=1}\mathbb{E}\left[\left(\int_Z
I_p(G_p^{q-1}f(z,\cdot)) \mu(dz)\right)^2\right],
$$
so that the estimate
$$
\sqrt{\mathbb{E}\left[\left(1-\frac{1}{q}\| DI_q(f)
\|^2_{L^2(\mu)}\right)^2\right]} \leq (\ref{bound_on_int})
$$
is proved, once we show that
\begin{equation}\label{bonI}
\mathbb{E}\left[\left(\int_Z I_p(G_p^{q-1}f(z,\cdot))
\mu(dz)\right)^2\right] \leq p! \int_{Z^p} \left\{\widehat{G}^q_p
f\right\}^2 d\mu^p.
\end{equation}
The proof of (\ref{bonI}) can be achieved by using the following
relations:
\begin{eqnarray}
&& \mathbb{E}\left[\left(\int_Z I_p(G_p^{q-1}f(z,\cdot))
\mu(dz)\right)^2\right] \notag\\
&&= \int_Z\int_Z\mathbb{E}[I_p(G_p^{q-1}f(z,\cdot)I_p(G_p^{q-1}f(z',\cdot)]\mu(dz)\mu(dz')\label{TMOE}\\
&& = p! \int_Z\int_Z \left[\int_{Z^p}
\widetilde{G_p^{q-1}f(z,\cdot)}\widetilde{G_p^{q-1}f(z',\cdot)}d\mu^p\right]\mu(dz)\mu(dz')\notag\\
&& = p!\int_{Z^p} \left[\int_Z
\widetilde{G_p^{q-1}f(z,\cdot)}\mu(dz)\right]^2 d\mu^p \leq
p!\int_{Z^p} \left[\int_Z G_p^{q-1}f(z,\cdot)\mu(dz)\right]^2
d\mu^p \notag\\
&& = p!\int_{Z^p} \left\{\widehat{G}^q_p f\right\}^2 d\mu^p.
\notag
\end{eqnarray}
Note that the use of the Fubini theorem in the equality
(\ref{TMOE}) is justified by the chain of inequalities
(\ref{Sco}), which is in turn a consequence of assumption
(\ref{technical}).

\noindent \underline{{\bf Step (S2)}}. First recall that
$$
\mathbb{E}\left[q^{-1}\|DI_q(f)\|^2_{L^2(\mu)}\right] =
\mathbb{E}[I_q(f)^2]=q!\|f\|_{L^2(\mu^q)}^2.
$$
Now use the Cauchy-Schwarz inequality, in order to write
\begin{eqnarray}
\frac{1}{q}\mathbb{E}\left[\int_Z |D_z I_q(f)|^3 \mu(dz)\right] &
\leq &
\frac{1}{q}\sqrt{\mathbb{E}\left[\|DI_q(f)\|^2_{L^2(\mu)}\right]}\times
\sqrt{\int_Z \mathbb{E}[(D_zI_q(f))^4]\mu(dz)} \label{D1}\\
&=& \sqrt{(q-1)!\|f\|_{L^2(\mu^q)}^2}\times \sqrt{\int_Z
\mathbb{E}[(D_zI_q(f))^4]\mu(dz)}. \label{D2}
\end{eqnarray}
By using (\ref{G0}) and (\ref{P0}), one deduces immediately that
$$
\{D_zI_q(f)\}^2 = q^2 \sum_{p=0}^{2(q-1)}
I_p[G_p^{q-1}f(z,\cdot)].
$$
As a consequence,
\begin{eqnarray} \label{S3_1}
\sqrt{\int_Z
\mathbb{E}[D_zI_q(f)^4]\mu(dz)}&\leq &q^2\sqrt{\sum_{p=0}^{2(q-1)} p!
\int_Z \left\{ \int_{Z^p} G_p^{q-1}f(z,\cdot)^2 d\mu^p \right\}
\mu(dz)} \\
& \leq & q^2\sum_{p=0}^{2(q-1)} p!^{1/2} \sqrt{\int_Z \left\{
\int_{Z^p} G_p^{q-1}f(z,\cdot)^2 d\mu^p \right\} \mu(dz)}, \notag
\end{eqnarray}
yielding the desired inequality.

\noindent \underline{\bf Step (S3)}. By using several times the
inequality $\sqrt{a+b} \leq \sqrt{a} + \sqrt{b}$ ($a,b\geq 0$) one
sees that, in order to prove
(\ref{starestimate1a})--(\ref{starestimate1b}), it is sufficient to
show that
\begin{eqnarray*}
&& \sum_{p=1}^{2(q-1)} p!^{1/2}\sqrt{ \int_{Z^p}
\left\{\widehat{G}^q_p
f\right\}^2 d\mu^p}  \\
&& \leq \sum_{t=1}^q \sum_{s=1}^{t\wedge(q-1)}\!\! \!{\bf 1}_{\{2
\leq t+s \leq 2q-1 \}}
(2q-t-s)!^{1/2}(t-1)!^{1/2}\binom{q-1}{t-1}\binom{t-1}{s-1}^{1/2}\!\!
\|f\star_t^s f \|_{L^2 (\mu^{2q-t-s})}.
\end{eqnarray*}
To see this, use (\ref{Gpihat3}) and the fact that (by
(\ref{SymInequality})) $\|\widetilde{f\star_t^s f} \|_{L^2
(\mu^{2q-t-s})} \leq \|f\star_t^s f \|_{L^2 (\mu^{2q-t-s})}$, to
obtain that
\begin{eqnarray*}
&& \sum_{p=1}^{2(q-1)} p!^{1/2}\sqrt{ \int_{Z^p}
\left\{\widehat{G}^q_p f\right\}^2 d\mu^p} \\
&& \leq \sum_{p=1}^{2(q-1)} p!^{1/2}\sum_{t=1}^{q}
\sum_{s=1}^{t\wedge (q-1) } {\bf 1} _{\{2q-t-s=p\}}
(t-1)!^{1/2}\binom{q-1}{t-1}\binom{t-1}{s-1} ^{1/2}\|f\star_t^s f
\|_{L^2 (\mu^{2q-t-s})},
\end{eqnarray*}
and then exploit the relation
$$
\sum_{p=1}^{2(q-1)} p!^{1/2}\sum_{t=1}^{q} \sum_{s=1}^{t\wedge (q-1)
} {\bf 1} _{\{2q-t-s=p\}} = \sum_{t=1}^q
\sum_{s=1}^{t\wedge(q-1)}\!\! \!{\bf 1}_{\{2 \leq t+s \leq 2q-1 \}}
(2q-t-s)!^{1/2} .
$$

\noindent\underline{{\bf Step (S4)}}. By using (\ref{oto}) and
some standard estimates, we deduce that
\begin{eqnarray}
&& \sum_{p=0}^{2(q-1)} p!^{1/2} \sqrt{\int_Z \left\{ \int_{Z^p}
G_p^{q-1}f(z,\cdot)^2 d\mu^p \right\} \mu(dz)} \notag \\
&& \leq \sum_{p=0}^{2(q-1)} p!^{1/2}\sum_{r=0}^{q-1}\sum_{l=0}^r
{\bf 1}_{\{ 2(q-1)-r-l= p \}} r!^{1/2}
\binom{q-1}{r}\binom{r}{l}^{1/2} \times  \label{Murray}\\
&& \text{\ \ \ \ \ \ \ \ \ \ \ \ \ \ \ \ \ \ \ \ \ \ \ \ \ \ \ \ \
\ \ \ \ \ \ \ \ \ \ \ \ \ \ \ }\times \sqrt{\int_Z \int_{Z^p}
[f(z,\cdot)\star^l_r f(z,\cdot)]^2 d\mu ^p \mu(dz)}.
\end{eqnarray}
We claim that, if (\ref{Declan}) is satisfied, then, for every
$r=0,...,q-1$ and $l=0,...,r$
\begin{equation}\label{keyEND}
\int_Z \int_{Z^p} [f(z,\cdot)\star^l_r f(z,\cdot)]^2 d\mu ^p
\mu(dz) =\int_{Z^{l+r+1}} [f\star_{q-l}^{q-1-r} f]^2 d\mu^{l+r+1}.
\end{equation}
In the two `easy' cases
\begin{itemize}
\item[(a)] $r=q-1$ and $l=1,...,q-1$,
\item[(b)] $r=1,...,q-1$ and $l=0$,
\end{itemize}
relation (\ref{keyEND}) can be deduced by a standard use of the
Fubini theorem (assumption (\ref{Declan}) is not needed here). Now
fix $p=1,...,2q-2$, as well as $r=1,...q-2$ and $l=0,...,r$ in
such a way that $2(q-1)-r-l=p$. For every fixed $z\in Z$, write
$|f(z,\cdot)|\star^l_r |f(z,\cdot)|$ to indicate the contraction
of indices $(r,l)$ obtained from the positive kernel
$|f(z,\cdot)|$. Note that, for $z$ fixed, such a contraction is a
function on $Z^p$, and also, in general, $|f(z,\cdot)|\star^l_r
|f(z,\cdot)|\geq f(z,\cdot)\star^l_r f(z,\cdot)$ and
$|f(z,\cdot)|\star^l_r |f(z,\cdot)|\neq f(z,\cdot)\star^l_r
f(z,\cdot)$. By a standard use of the Cauchy-Schwarz inequality
and of the Fubini theorem, one sees that
\begin{equation}\label{cr}
\int_Z \int_{Z^p} \left[|f(z,\cdot)|\star^l_r
|f(z,\cdot)|\right]^2 d\mu ^p \mu(dz)\leq \|f\star_q^{q-1-r+l} f
\|^2_{L^2(\mu^{r-l-1})}< \infty,
\end{equation}
where the last relation is a consequence of assumption
(\ref{Declan}) as well as of the fact that, by construction,
$1\leq 1+r-l\leq q-1$. Relation (\ref{cr}) implies that one can
apply the Fubini Theorem to the quantity
$$
\int_Z \int_{Z^p} [f(z,\cdot)\star^l_r f(z,\cdot)]^2 d\mu ^p
\mu(dz)
$$
(by first writing the contractions $ f(z,\cdot)\star^l_r
f(z,\cdot)$ in an explicit form), so to obtain the desired
equality (\ref{keyEND}). By plugging (\ref{keyEND}) into
(\ref{Murray}) and by applying the change of variables $a=q-1-r$
and $b=q-l$, one deduces
(\ref{starestimate2a})--(\ref{starestimate2b}). This concludes the
proof of Theorem \ref{T : MWII-Bounds}.
 \fin

\section{Central limit theorems}\label{SS : CLTchaos}\setcounter{equation}{0}
We consider here CLTs. The results of
this section generalize the main findings of \cite{PecTaq07}.
The following result uses Theorem \ref{T : MWII-Bounds} in order
to establish a general CLT for multiple integrals of arbitrary
order.

\begin{thm}[CLTs on a fixed chaos]\label{T : chaosCLT} Let $X\sim \mathscr{N}(0,1)$.
Suppose that $\mu(Z)=\infty$, fix $q\geq 2$, and let $F_k =
I_q(f_k)$, $k\geq 1$, be a sequence of multiple stochastic
Wiener-Itô integrals of order $q$. Suppose that, as
$k\rightarrow\infty$, the normalization condition
$\mathbb{E}(F_k^2) = q!\|f_k\|^2_{L^2(\mu^q)} \rightarrow 1$ takes
place. Assume moreover that the following three conditions hold:
\begin{itemize}
\item[\rm (I)] For every $k\geq 1$, the kernel $f_k$ verifies
(\ref{technical}) for every $p=1,...,2(q-1)$.
%\item[\rm (II)] For every $k\geq 1$, every $z\in Z$, every $r=1,...,q-1$ and
%every $l=0,...,r$ such that $r\neq l$, the kernel
%$f_k(z,\cdot)\star_{ r}^{ l} f_k(z,\cdot)$ is an element of
%$L^2_s(\mu^{2(q-1)-r-l})$.
\item[\rm (II)] For every $r=1,...,q$, and every $l=1,...,r \!
\wedge \! (q-1)$, one has that $f_k \star_r^l f_k \! \in \!
L^2(\mu^{2q-r-l})$ and also $\|f_k \star_r^l f_k
\|_{L^2(\mu^{2q-r-l})} \rightarrow 0$ (as $k\rightarrow \infty$).
\item[\rm (III)] For every $k\geq 1$, one has that $\int_{Z^q}
f^4_k\,\, d\mu ^q <\infty$ and, as $k\rightarrow\infty$,
$\int_{Z^q} f^4_k\,\, d\mu ^q \rightarrow 0$.
\end{itemize}
Then, $F_k \stackrel{\rm law}{\rightarrow} X$, as $k\rightarrow
\infty$, and formulae (\ref{starestimate1a})--(\ref{revolved3})
provide explicit bounds in the Wasserstein distance
$d_W(I_q(f_k),X)$.
\end{thm}
{\bf Proof.} First note that the fact that $f_k \star_r^l f_k \in
L^2(\mu^{2q-r-l})$ for every $r=1,...,q$ and every $l=1,...,r\wedge
(q-1)$ (due to Assumption (II)) imply that Assumption (ii) in the
statement of Theorem \ref{T : MWII-Bounds} holds for every $k\geq 1$
(with $f_k$ replacing $f$), and also (by using (\ref{revolved2}))
that condition (\ref{Declan}) is satisfied by each kernel $f_k$. Now
observe that, if Assumptions (I)-(III) are in order and if
$\mathbb{E}(F_k^2) = q!\|f_k\|^2_{L^2(\mu^q)} \rightarrow 1$, then
relations (\ref{starestimate1a})--(\ref{revolved3}) imply that
$d_W(F_k,X)\rightarrow 0$. Since convergence in the Wasserstein
distance implies convergence in law, the conclusion is immediately
deduced. \fin

\begin{example}\label{EX : MWII_order2_CLT}
{\rm Consider a sequence of double integrals of the type $I_2(f_k)$, $k\geq 1$
where $f_k \in L_s^2(\mu^2)$ satisfies (\ref{technical_q=2}) (according
to Remark \ref{R : Tech}(1), this implies that (\ref{technical}) is
satisfied). We suppose that the following three conditions are
satisfied: (a) $\mathbb{E}I_2(f_k)^2 = 2\|f_k\|_{L^2(\mu^2)}^2 = 1$, (b)
$f_k\star_2^1 f_k\in L^2 (\mu^{1})$ and (c) $f_k\in L^4(\mu^2)$. Then, according to Theorem \ref{T : chaosCLT}, a sufficient condition in order to have that (as $k\rightarrow \infty$)
\begin{equation}\label{MR0}
I_2(f_k) \stackrel{\rm law}{\rightarrow} \mathscr{N}(0,1)
\end{equation}
is that
\begin{eqnarray}\label{MR}
&& \|f_k\|_{L^4(\mu^2)} \rightarrow 0, \\ \label{MR2}
&& \|f_k\star_2^1 f_k\|_{L^2(\mu)} \rightarrow 0,\, \, {\rm and} \,\, \|f_k\star_1^1 f_k\|_{L^2(\mu^2)}\rightarrow 0.
 \end{eqnarray}
This last fact coincides with the content of Part 1 of Theorem 2 in \cite{PecTaq07}, where one can find an alternate proof based on a decoupling technique, known as the ``principle of conditioning'' (see e.g. Xue \cite{XUE}). Note that an explicit upper bound for the Wasserstein distance can be deduced from relation (\ref{doubleEX}).} \fin
\end{example}

\begin{example}\label{E : MWII-3-p2} {\rm Consider a sequence of random
variables of the type $I_3(f_k)$, $k\geq 1$, with unitary variance
and verifying Assumption (I) in the statement of Theorem \ref{T :
chaosCLT}. Then, according to the conclusion of Theorem \ref{T :
chaosCLT}, a sufficient condition in order to have that, as
$k\rightarrow \infty$,
$$I_3(f_k)\stackrel{\rm law}{\longrightarrow} X
\sim\mathscr{N}(0,1),$$ is that the following  six quantities
converge to zero:
\begin{eqnarray*}
&& \|f_k\star_1^1 f_k\|_{L^2(\mu^4)}, \,\,\, \|f_k\star_2^1
f_k\|_{L^2(\mu^3)}, \,\,\, \|f_k\star_2^2 f_k\|_{L^2(\mu^2)}, \\
&& \|f_k\star_3^1  f_k \|_{L^2(\mu^2)}, \,\,\, \|f_k\star_3^2 f_k
\|_{L^2(\mu^1)},\, \text{  \ and \ }\, \|f_k\|_{L^4(\mu^3)}.
\end{eqnarray*}
Moreover, an explicit upper bound in the Wasserstein distance is
given by the estimate (\ref{ZZ}).} \fin
\end{example}
%\subsection{Necessary and sufficient conditions in the case of double integrals} \label{SS : WeakSingle2ble}
The following result, proved in \cite[Theorem 2]{PecTaq07}, represents a counterpart to the CLTs for double integrals discussed in Example \ref{EX : MWII_order2_CLT}.
\begin{prop} [See \cite{PecTaq07}]
\label{T : PoissCLT} Consider a sequence $F_k = I_{2}(f_k)$, $k\geq 1$, of double integrals verifying assumptions {\rm (a), (b)} and {\rm (c)} of Example \ref{EX : MWII_order2_CLT}. Suppose moreover that (\ref{MR}) takes place. Then,
\begin{enumerate}
\item[$1.$] if $F_{k}\in L^{4}\left( \mathbb{P}\right) $ for every $k$,
a sufficient condition to have (\ref{MR2}) is that%
\begin{equation}
\mathbb{E}\left( F_k^{4}\right) \rightarrow 3;  \label{GGG}
\end{equation}
\item[$2.$] if the sequence $\left\{F_{k}^{4}:k \geq 1\right\} $ is
uniformly integrable, then conditions (\ref{MR0}), (\ref{MR2}) and
(\ref{GGG}) are equivalent.
\end{enumerate}
\end{prop}

For the time being, it seems quite hard to prove a result analogous to Proposition \ref{T : PoissCLT} for a sequence of multiple integrals of order $ q \geq 3$.

\section{Sum of a single and a double integral}\label{SS : SDint}
\setcounter{equation}{0}
As we will see in the forthcoming Section \ref{S : EX}, when dealing with
quadratic functionals of stochastic processes built from completely
random measures, one needs explicit bounds for random variables of
the type $F = I_1(g) + I_2 (h)$, that is, random variables that are
the sum of a single and a double integral. The following result,
that can be seen as a generalization of Part B of Theorem 3 in \cite{PecTaq07}, provides explicit bounds for random variables of this
type.

\begin{thm}\label{T : Bound1C2C}
Let $F = I_1(g) + I_2(h)$ be such that
\begin{itemize}
\item[\rm (I)] The function $g$ belongs to $L^2(\mu) \cap
L^3(\mu)$;
\item[\rm (II)] The kernel $h\in L^2_s (\mu^2)$ is such that: (a) $h\star_2^1 h\in L^2 (\mu^{1})$, (b) relation (\ref{technical_q=2}) is verified, with $h$ replacing $f$, and (c) $h\in L^4(\mu^2)$.
\end{itemize}
Then, one has the following upper bound on the Wasserstein
distance between the law of $F$ and the law of
$X\sim\mathscr{N}(0,1)$:
\begin{eqnarray}\label{1+2Bounds}
&&d_W(F,X) \\
&& \leq  \left |1 - \|g\|^2_{L^2(\mu)}-2\|h\|^2_{L^2(\mu^2)}
\right|
+2\| h\star_2^1 h \|_{L^2(\mu)} \notag \\
&& + \sqrt{8}\| h\star_1^1 h \|_{L^2(\mu^2)} + 3 \|g\star_1^1 h
\|_{L^2(\mu)} + 32 \|g\|^3_{L^3(\mu)} \notag \\
&& + 4 \|h\|_{L^2(\mu^2)} \times
\left\{ \|h\|^2_{L^4(\mu^2)} + 2^{1/2} \|h\star_2^1 h\|_{L^2(\mu)}
\right\}.\notag
\end{eqnarray}
The following inequality also holds:
\begin{equation}\label{sub}
\|g\star_1^1 h \|_{L^2(\mu)} \leq \|g\|_{L^2(\mu)} \times
\|h\star_1^1 h\|^{1/2}_{L^2(\mu^2)}.
\end{equation}
\end{thm}
{\bf Proof.} Thanks to Theorem \ref{T : MainUpperBound}, we know
that $d_W (F , X)$ is less or equal to the RHS of
(\ref{GenUpBound}). We also know that $$ D_z F = g(z)
+2I_1(h(z,\cdot)), \text{ \ and \ } -D_z L^{-1} F = g(z)
+I_1(h(z,\cdot)).$$ By using the multiplication formula
(\ref{pproduct}) (in the case $p=q=1$) as well as a Fubini argument,
one easily deduces that
$$ \int_Z D_zF \times (-D_z L^{-1}F) \mu(dz) =
\|g\|^2_{L^2(\mu)}+2\|h\|^2_{L^2(\mu^2)} \! +\! 2I_1(h\star_2^1 h) +
2I_2(h\star_1^1 h)\!+\! 3I_1 (g\star_1^1 h). $$ This last relation
yields
\begin{eqnarray}\label{OLL} && \sqrt{\mathbb{E}\left[(1-\langle DF,
-DL^{-1}F\rangle _{L^2(\mu)})^2\right]} \\ &&\leq   \left |1 -
\|g\|^2_{L^2(\mu)}-2\|h\|^2_{L^2(\mu^2)} \right|  +2\| h\star_2^1 h
\|_{L^2(\mu)}  + \sqrt{8}\| h\star_1^1 h \|_{L^2(\mu^2)} + 3
\|g\star_1^1 h \|_{L^2(\mu)}. \notag
\end{eqnarray}
To conclude the proof, one shall use the following relations,
holding for every real $a,b$:
\begin{equation}\label{chain} (a + 2b)^2|a+b| \leq ( |a| + 2|b| )^2 (|a|+|b|) \leq ( |a| + 2|b| )^3 \leq 4|a|^3 + 32|b|^3.
\end{equation}
By applying (\ref{chain}) in the case $a = g(z)$ and
$b=I_1(h(z,\cdot))$, one deduces that
\begin{eqnarray}\label{ss}
&&\int_Z \mathbb{E}\left[|D_z F|^2 |D_z L^{-1}F|\right]\mu(dz)  \\
&& \leq 4 \mathbb{E}\int_Z |I_1(h(z,\cdot))|^3\mu(dz) + 32\int_Z
|g(z)|^3\mu(dz). \notag
\end{eqnarray}
By using the Cauchy-Schwarz inequality, one infers that
\begin{eqnarray} \notag
&& 4 \mathbb{E}\int_Z
|I_1(h(z,\cdot))|^3\mu(dz)\\
&& \leq \sqrt{\mathbb{E}\int_Z |I_1(h(z,\cdot))|^4 \mu(dz)} \times
4\|h\|_{L^2(\mu^2)},
\end{eqnarray}
and (\ref{1+2Bounds}) is deduced from the equality
$$
\mathbb{E}\int_Z |I_1(h(z,\cdot))|^4 \mu(dz) = 2\|h\star_2^1
h\|^2 _{L^2(\mu)} +\|h\|^4_{L^4(\mu^2)}.
$$
Formula (\ref{sub}) is once again an elementary consequence of the
Cauchy-Schwarz inequality. \fin

%\section{Applications to functionals of Lévy processes} \label{Levy}
%
%{\bf *** HERE, I think we should try to translate most of the
%results in terms of Lévy processes, in particular so to establish
%an exhaustive connection with the paper \cite{SolUtzViv}. It seems
%to me that this is would enhance very much the applicability of
%our results, even outside the scope of the paper! WHat do you
%think? *** }

\begin{rem}{\rm Consider a sequence of vectors $(F_k, H_k)$, $k\geq 1$ such that: (i) $F_k = I_2(f_k)$, $k\geq 1$, is a sequence of double integrals verifying assumptions (a), (b) and (c) in Example \ref{EX : MWII_order2_CLT}, and (ii) $H_k=I_1(h_k)$, $k\geq 1$, is a sequence of single integrals with unitary variance. Suppose moreover that the asymptotic relations in (\ref{1CONDclt}), (\ref{MR}) and (\ref{MR2}) take place. Then, the estimates
(\ref{1+2Bounds})--(\ref{sub}) yield that, for every $ (\alpha , \beta ) \neq (0,0)$, the Wasserstein distance
between the law of
$$\frac{1}{\sqrt{\alpha^2 + \beta^2}}\left( \alpha  F_k + \beta
H_k\right)
$$
and the law of $X \sim\mathscr{N}(0,1)$, converges to zero as
$k\rightarrow \infty$, thus implying that $(F_k, H_k)$ converges in law to a vector $(X,X')$ of i.i.d. standard Gaussian random variables. Roughly speaking, this last fact implies that, when assessing the asymptotic joint Gaussianity of a vector such as $(F_k,H_k)$, one can study \textsl{separately} the one-dimensional sequences $\{F_k\}$ and $\{H_k\}$. This phenomenon coincides with the content of Part B of Theorem 3 in \cite{PecTaq07}. See \cite{NouPeRev}, \cite{NO} and \cite{PT} for similar results involving vectors of multiple integrals (of arbitrary order) with respect to Gaussian random measures.}
\end{rem}

\section{Applications to non-linear functionals of Ornstein-Uhlenbeck Lévy processes} \label{S : EX}
As an illustration, in this section we focus on CLTs related to Ornstein-Uhlenbeck Lévy processes, that is, processes obtained by integrating an exponential kernel of the type $$x \rightarrow \sqrt{2\lambda}\times{\rm e} ^{ -\lambda(t-x)}{\bf 1}_{\{x\leq t\}},$$ with respect to an independently scattered random measure. Ornstein-Uhlenbeck Lévy processes have been recently applied to a variety of frameworks, such as finance (where they are used to model stochastic volatility -- see e.g. \cite{BNS}) or non-parametric Bayesian survival analysis (where they represent random hazard rates -- see e.g. \cite{DBPP, nieto05, PePrU}). In particular, in the references \cite{DBPP} and \cite{PePrU} it is shown that one can use some of the CLTs of this section in the context of Bayesian prior specification.

\subsection{Quadratic functionals of Ornstein-Uhlenbeck processes}\label{SS : QuadOU}
\setcounter{equation}{0}
We consider the stationary \textsl{\ Ornstein-Uhlenbeck L\'{e}vy
process} given by
\begin{equation}
Y_{t}^{\lambda }=\sqrt{2\lambda }\int_{-\infty
}^{t}\int_{\mathbb{R}}u\exp \left( -\lambda \left( t-x\right)
\right) \widehat{N}\left( du,dx\right) \text{, \ \ }t\geq 0\text{,}
\label{OUl}
\end{equation}%
where $\widehat{N}$ is a centered Poisson measure over
$\mathbb{R\times R}$, with control measure given by $\nu \left(
du\right) dx$, where $\nu \left( \cdot \right) $ is positive, non-atomic and $\sigma$-finite. We
assume also that $\int u^{j}\nu \left( du\right) <\infty $ for
$j=2,3,4,6$, and $\int u^{2}\nu \left( du\right) =1$. In particular,
these assumptions yield that
\begin{equation*}
\mathbb{E}\left[\left( Y_{t}^{\lambda }\right) ^{2}\right]=2\lambda
\int_{-\infty }^{t}\int_{\mathbb{R}}u^{2}e^{-2\lambda \left(
t-x\right) }\nu \left( du\right) dx=1.
\end{equation*}
The following result has been proved in \cite[Theorem 5]{PecTaq07}.

\begin{thm}[See \cite{PecTaq07}]
\label{P : OUclt} For every $\lambda >0$, as $T\rightarrow \infty$,
\begin{equation}
Q(T,\lambda) := \sqrt{T}\left\{ \frac{1}{T}\int_{0}^{T}\left(
Y_{t}^{\lambda }\right) ^{2}dt-1\right\} \overset{\rm
law}{\longrightarrow }\sqrt{\frac{1}{\lambda }+c_{\nu }^{2}}\times
X\text{,} \label{OU_CLT!}
\end{equation}%
where $c^2_\nu  := \int u^4 \nu(du)$ and $X\sim \mathscr{N}\left( 0,1\right) $ is a centered standard
Gaussian random variable.
\end{thm}

By using Theorem \ref{T : Bound1C2C}, one can obtain the following
Berry-Ess\'{e}en estimate on the CLT appearing in (\ref{OU_CLT!})
(compare also with Example \ref{EX : OUlinear}).

\begin{thm}\label{T : UpperBouhl} Let $Q(T,\lambda)$, $T>0$, be
defined as in (\ref{OU_CLT!}), and set $$\tilde{Q}(T,\lambda) :=
Q(T,\lambda)/\sqrt{\frac{1}{\lambda }+c_{\nu }^{2}}.$$ Then, there
exists a constant $0<\gamma(\lambda)<\infty$, independent of $T$ and
such that
\begin{equation}
d_W(\tilde{Q}(T,\lambda),X) \leq \frac{\gamma(\lambda)}{\sqrt{T}}.
\end{equation}
\end{thm}

{\bf Proof}. Start by introducing the notation
\begin{eqnarray}
H_{\lambda ,T}\left( u,x;u,^{\prime }x^{\prime }\right) &=&\left(
u\times u^{\prime }\right) \frac{\mathbf{1}_{\left( -\infty
,T\right] ^{2}}\left( x,x^{\prime }\right) }{T\sqrt{\frac{1}{\lambda
}+c_{\nu }^{2}}}\left\{ e^{\lambda \left( x+x^{\prime }\right)
}\left( 1-e^{-2T}\right) \mathbf{1}_{\left( x\vee x^{\prime }\leq
0\right)
}+\right.  \label{def} \\
&&\text{ \ \ \ \ \ \ \ \ \ \ \ \ \ \ \ \ \ \ \ \ \ \ \ \ \ }+\left.
e^{\lambda \left( x+x^{\prime }\right) }\left( e^{-2\lambda \left(
x\vee x^{\prime }\right) }-e^{-2\lambda T}\right) \mathbf{1}_{\left(
x\vee
x^{\prime }>0\right) }\right\} \text{,}  \notag \\
H_{\lambda ,T}^{\ast }\left( u,x\right)
&=&u^{2}\frac{\mathbf{1}_{\left( -\infty ,T\right] }\left( x\right)
}{T\sqrt{\frac{1}{\lambda }+c_{\nu }^{2}}}\left\{ e^{2\lambda
x}\left( 1-e^{-2T}\right) \mathbf{1}_{\left( x\leq 0\right)
}+e^{2\lambda x}\left( e^{-2\lambda x}-e^{-2\lambda T}\right)
\mathbf{1}_{\left( x>0\right) }\right\} \text{.}  \notag
\end{eqnarray}%
As a
consequence of the multiplication formula (\ref{pproduct}) and of a
standard Fubini argument, one has (see \cite[Proof of Theorem 5]{PecTaq07})
$$
\tilde{Q}(T,\lambda) = I_1\left(\sqrt{T}H_{\lambda ,T}^{\ast
}\right) + I_2\left(\sqrt{T}H_{\lambda ,T} \right),
$$
a combination of a single and of a double integral. To apply Theorem \ref{T : Bound1C2C},
we use the following asymptotic
relations (that one can verify by resorting to the explicit
definitions of $H_{\lambda ,T}$ and $ H^{\ast}_{\lambda ,T}$ given
in (\ref{def})), holding for $T\rightarrow \infty$:
\begin{itemize}
\item[(a)] $$\left|1 -\left \|\sqrt{T}H_{\lambda ,T}^{\ast
}\right\|^2_{L^2(d\nu dx)}- 2\left\|\sqrt{T}H_{\lambda
,T}\right\|^2_{L^2((d\nu dx)^2)} \right|= O\left(\frac1T\right)$$
\item[(b)] $$ \left\|\sqrt{T}H_{\lambda ,T}^{\ast
}\right \|^3_{L^3(d\nu dx)} \sim \frac{1}{\sqrt{T}} $$
\item[(c)] $$ \left\|\sqrt{T}H_{\lambda ,T}\right \|^2_{L^4((d\nu dx))^2} \sim \frac{1}{\sqrt{T}}$$
\item[(d)] $$ \left\|(\sqrt{T}H_{\lambda ,T})\star_2^1 (\sqrt{T}H_{\lambda
,T})
\right \|_{L^2(d\nu dx)} \sim \frac{1}{\sqrt{T}}$$
\item[(e)] $$ \left\|(\sqrt{T}H_{\lambda ,T})\star_1^1(\sqrt{T}H_{\lambda ,T}
) \right \|_{L^2((d\nu dx)^2)} \sim \frac{1}{\sqrt{T}}$$
\item[(f)] $$ \left\|(\sqrt{T}H^{\ast}_{\lambda ,T})\star_1^1(\sqrt{T}H_{\lambda ,T}
) \right \|_{L^2(d\nu dx)} \sim \frac{1}{\sqrt{T}}.$$
\end{itemize}
The conclusion is obtained by using the estimates
(\ref{1+2Bounds})--(\ref{sub}) and applying Theorem \ref{T : Bound1C2C}.
\fin
\subsection{Berry-Ess\'{e}en bounds for arbitrary tensor powers of Ornstein-Uhlenbeck kernels}\label{SS : tensorOU}

Let $\widehat{N}$ be a centered Poisson measure over
$\mathbb{R\times R}$, with control measure given by $\mu(du,dx) = \nu \left(
du\right) dx$, where $\nu \left( \cdot \right) $ is positive, non-atomic and $\sigma$-finite. We
assume that $\int u^{2}\nu \left( du\right) =1$ and $\int u^{4}\nu \left( du\right) <\infty $. Fix $\lambda >0$, and, for every $t\geq 0$, define the Ornstein-Uhlenbeck kernel
\begin{equation}\label{ft}
f_t(u,x) = u\times\sqrt{2\lambda} \exp\{-\lambda(t-x)\} {\bf 1}_{\{x\leq t\}}, \,\,\, (u,x)\in \R \times\R.
\end{equation}
For every fixed $q\geq 2$, we define the $q$th {\sl tensor power} of $f_t$, denoted by $f_t^{\otimes q}$, as the symmetric kernel on $(\R \times\R)^q$ given by

\begin{equation}\label{ftq}
f_t^{\otimes q} (u_1,x_1;...;u_q,x_q)=\prod_{j=1}^q f_t(u_j,x_j).
\end{equation}

We sometimes set $y=(u,x)$. Note that, for every $t\geq 0$, one has that $\int f^2_t(y) \mu(dy) = \int f^2_t(u,x)\nu(du)dx = 1$, and therefore $f_t^{\otimes q} \in L_s^2(\mu^q)$; it follows that the multiple integral $$Z_t(q): = I_q (f_t^{\otimes q})$$ is well defined for every $t\geq 0$.

\begin{rem}
{\rm Fix $q\geq 2$, and suppose that $\int |u|^{j} d\nu <\infty$, $\forall j=1,...,2q$. Then, one can prove that $Z_t(1) ^q$ is square-integrable and also that the random variable $Z_t(q)$ coincides with the projection of $Z_t(1) ^q $ on the $q$th Wiener chaos associated with $\widehat{N}$. This fact can be easily checked when $q=2$: indeed (using Proposition \ref{P : ProductPoisson} in the case $f=g=f_t$) one has that $$ Z_t(1)^2 = I_1(f_t^2) + I_2(f_t^{\otimes 2}) = I_1(f_t^2) +Z_t (2) ,$$ thus implying the desired relation. The general case can be proved by induction on $q$.}
\end{rem}

The main result of this section is the following application of Theorem \ref{T : MWII-Bounds} and Theorem \ref{T : chaosCLT}.
\begin{thm}\label{T : TensOU} Fix $\lambda >0$ and $q\geq 2$, and define the positive constant $c = c(q,\lambda):= 2(q-1)!/\lambda$. Then, one has that, as $T\rightarrow \infty$,
\begin{equation}\label{CLTout}
M_T(q) := \frac1{\sqrt{cT}}\int_0^T Z_t (q) dt \stackrel{\rm law}{\longrightarrow} X\sim\mathscr{N} (0,1),
\end{equation}
and there exists a finite constant $\rho = \rho(\lambda, q, \nu)>0$ such that, for every $T>0$,
\begin{equation}\label{BEout}
d_W(M_T(q),\, X) \leq \frac\rho{\sqrt{T}}.
\end{equation}
\end{thm}
{\bf Proof.} The crucial fact is that, for each $T$, the random variable $M_T(q)$ has the form of a multiple integral, that is, $M_T(q) = I_q (F_T)$, where $F_T \in L_s^2 (\mu ^q)$ is given by $$F_T(u_1,x_1;...;u_q,x_q) =\frac{1}{\sqrt{cT} } \int _0 ^T f_t^{\otimes q} (u_1,x_1;...;u_q,x_q)dt,$$ where $f_t^{\otimes q}$ has been defined in (\ref{ftq}). By using the fact that the support of $F_T$ is contained in the set $(\R \times (-\infty,T] )^q$ as well as the assumptions on the second and fourth moments of $\nu$, one easily deduces that the technical condition (\ref{technical}) (with $F_T$ replacing $f$) is satisfied for every $T\geq 0$. According to Theorem \ref{T : MWII-Bounds} and Theorem \ref{T : chaosCLT}, both claims (\ref{CLTout}) and (\ref{BEout}) are proved, once we show that, as $T\rightarrow\infty$, one has that
\begin{equation}\label{outEXP}
|1 -\mathbb{E}(M_T(q))^2| \sim 1/T,
\end{equation}
and also that
\begin{eqnarray}\label{out4th}
&& \|F_T\|^2_{L^4(\mu^q)} = O(1/T), \,\, {\rm and} \\ \label{outSTAR}
&& \|F_T \star _r ^l F_T \|_{L^2(\mu^{2q-r-l})} = O(1/T), \,\, \forall r=1,...,q,\,\, \forall l=1,...,r\wedge(q-1)
\end{eqnarray}
(relation (\ref{outEXP}) improves the bounds in Theorem \ref{T : MWII-Bounds}). In order to prove (\ref{outEXP})--(\ref{outSTAR}), for every $t_1,t_2\geq 0$ we introduce the notation
\begin{equation}\label{angle}
\langle f_{t_1} , f_{t_2} \rangle_\mu = \int_{\R\times\R} f_{t_1}(y) f_{t_2}(y)\, \mu(dy) = {\rm e}^{-\lambda(t_1 +t_2)} {\rm e}^{2 \lambda(t_1 \wedge t_2)},
\end{equation}
(recall that $\int u^2 d\nu =1$) and also, for $t_1,t_2,t_3,t_4 \geq 0$,
\begin{eqnarray}\label{angle1}
\langle f_{t_1} f_{t_2}\, , \,f_{t_3} f_{t_4} \rangle_\mu &=& \int_{\R\times\R} f_{t_1}(y)f_{t_2}(y)f_{t_3}(y)f_{t_4} (y)\,\, \mu(dy) \\ \label{angle2} &=&\left[\int_\R u^4 \nu(du)\right] \times \lambda {\rm e}^{-\lambda(t_1 +t_2+t_3+t_4)} {\rm e}^{4 \lambda(t_1 \wedge t_2\wedge t_3 \wedge t_4)}.
\end{eqnarray}
To prove (\ref{outEXP}), one uses the relation (\ref{angle}) to get
$$
\mathbb{E}\left[M_T(q)^2\right] = \frac {q!}{cT} \int_0^T\int_0^T \langle f_{t_1}, f_{t_2}\rangle_\mu ^q \, dt_1 dt_2 = 1- \frac{1}{Tq\lambda} (1-{\rm e}^{-\lambda q T}).
$$
In the remaining of the proof, we will write $\kappa$ in order to indicate a strictly positive finite constant independent of $T$, that may change from line to line. To prove (\ref{out4th}), one uses the fact that
$$
\int_{(\R \times \R)^q} F_T^4 d\mu^q = \frac{1}{c^2 T^2} \int_0^T\int_0^T\int_0^T\int_0^T \langle f_{t_1} f_{t_2}\, , \,f_{t_3} f_{t_4} \rangle_\mu ^q \, \, dt_1 dt_2 dt_3 dt_4 \leq \frac{\kappa}{T},
$$
where the last relation is obtained by resorting to the explicit representation (\ref{angle2}), and then by evaluating the restriction of the quadruple integral to each simplex of the type $\{t_{\pi(1)} >t_{\pi(2)} >t_{\pi(3)} > t_{\pi(4)}\}$, where $\pi$ is a permutation of the set $\{1,2,3,4\}$. We shall now verify the class of asymptotic relations (\ref{outSTAR}) for $r=q$ and $l=1,...,q$. With $y=(u,x)$, one has
$$ F_T\star^l_q F_T (y_1,...,y_{q-l}) = \frac{1}{cT}\int_0^T\int_0^T \left[\prod_{i=1}^{q-l} f_{t_1}(y_i)f_{t_2}(y_i)\right]\langle f_{t_1}f_{t_2}\rangle ^l_\mu \,\, dt_1 dt_2,$$
and hence
\begin{eqnarray*}
&& \|F_T \star _q ^l F_T \|_{L^2(\mu^{q-l})}^2 \\
&& = \frac{\kappa}{T^2} \int_0^T\int_0^T\int_0^T\int_0^T \langle f_{t_1} f_{t_2}\, , \,f_{t_3} f_{t_4} \rangle_\mu ^{q-l} \langle f_{t_1} f_{t_2} \rangle_\mu^l \langle f_{t_3} f_{t_4} \rangle_\mu^l \, \, dt_1 dt_2 dt_3 dt_4\\
&& \leq \frac{\kappa}{T},
\end{eqnarray*}
where the last relation is verified by first using (\ref{angle})--(\ref{angle2}), and then by assessing the restriction of the quadruple integral to each one of the $4!=24$ simplexes of the type $\{t_{\pi(1)} >t_{\pi(2)} >t_{\pi(3)} > t_{\pi(4)}\}$. To deal with (\ref{outSTAR}) in the case $r=1,...,q-1$ and $l=1,...,r$, one uses the fact that
\begin{eqnarray*}
&&F_T\star^l_r F_T (y_1,...,y_{r-l},w_1,...,w_{q-r},z_1,...,z_{q-r}) \\
&& = \frac{1}{cT}\int_0^T\int_0^T \left[\prod_{i=1}^{r-l} f_{t_1}(y_i)f_{t_2}(y_i)\right]\times \\
&& \quad\quad\quad\quad\quad\quad\quad\quad\quad \times \left(f_{t_1}(w_1)\cdot\cdot\cdot f_{t_1}(w_{q-r})f_{t_2}(z_1)\cdot\cdot\cdot f_{t_2}(z_{q-r})\right)\langle f_{t_1}f_{t_2}\rangle ^l_\mu \,\, dt_1 dt_2,
\end{eqnarray*}
and therefore
\begin{eqnarray*}
&& \|F_T \star _r ^l F_T \|_{L^2(\mu^{2q-r-l})}^2 \\
&& = \frac{\kappa}{T^2} \int_0^T\!\int_0^T\!\int_0^T\!\int_0^T\! \langle f_{t_1} f_{t_2}\, , \,f_{t_3} f_{t_4} \rangle_\mu ^{r-l} \langle f_{t_1} f_{t_3} \rangle_\mu^{q-r} \langle f_{t_2} f_{t_4} \rangle_\mu^{q-r} \langle f_{t_1} f_{t_2} \rangle_\mu^l \langle f_{t_3} f_{t_4} \rangle_\mu^l \, \, dt_1 dt_2 dt_3 dt_4\\
&& \leq \frac{\kappa}{T},
\end{eqnarray*}
where the last relation is once again obtained by separately evaluating each restriction of the quadruple integral over a given simplex. This concludes the proof.
\fin
{\bf Acknowledgments.} Part of this paper has been written while G. Peccati was visiting the Department of Mathematics of the ``Universitat Aut\'{o}noma de Barcelona'', in April 2008. This author heartily thanks J.L. Sol\'{e} and F. Utzet for their hospitality and support. J.L. Sol\'{e} and F. Utzet were supported by the Grant BFM2006-06427, Ministerio de Educaci\'{o}n y Ciencia and FEDER. M.S. Taqqu was partially supported by the NSF Grants DMS-0505747 and DMS-0706786 at Boston University.


\begin{thebibliography}{99}

%\bibitem{BGK} Bhansalia, R.J., Giraitis, L. and Kokoszka, P.S. (2007).
%Approximations and limit theory for quadratic forms of linear processes.
%\textit{Stochastic Processes and their Applications}, \textbf{117}, 71-95.
%
%\bibitem{BGK2} Bhansalia, R.J., Giraitis, L. and Kokoszka, P.S. (2007).
%Convergence of quadratic forms with nonvanishing diagonals. \textit{%
%Statistics and Probability Letters}, to appear.

\bibitem{BNS} Bandorff-Nielsen, O.E. and Shepard, N. (2001). Non Gaussian
Ornstein-Uhlenbeck-based models and some of their uses in financial
economics. \textit{Journal of the Royal Statistical Society B} \textbf{63}, 167--241.

\bibitem{Chatterjee_AOP} Chatterjee, S. (2007). A new method of normal approximation.
To appear in: \textit{The Annals of Probability}.
%
%\bibitem{Chatterjee_ptrf} Chatterjee, S. (2007). Fluctuation of eigenvalues and second order Poincaré inequalities.
%To appear in: \textit{Probability Theory and Related Fields}.

\bibitem{Chen_Shao_sur} Chen, L. and Shao, Q.-M.  (2005).
Stein's method for normal approximation. In: \textit{An
introduction to Stein's method}, 1-59. Lect. Notes Ser. Inst.
Math. Sci. Natl. Univ. Singap. \textbf{4}, Singapore Univ. Press,
Singapore, 2005.

%\bibitem{CohTaq} Cohen, S. and Taqqu, M.S. (2004). Small and large scale
%behavior of the Poissonized Telecom process. \textit{Methodology and
%Computing in Applied Probability }\textbf{6}, 363-379.

\bibitem{DBPP} De Blasi, P., Peccati, G. and Pr\"{u}nster, I. (2008).
Asymptotics for posterior hazards. To appear in: \textit{The
Annals of Statistics.}

\bibitem{Dec_Savy} Decreusefond, L. and Savy, N. (2007). Rubinstein
distance on configuration space. Preprint.

%\bibitem{DeJong} De Jong, P. (1987) A central limit theorem for generalized
%quadratic forms. \textit{Probability Theory and Related Fields}, \textbf{75}%
%, 261-277.

\bibitem{Dudely} Dudley, R.M. (2002). \textit{Real Analysis and Probability.}
Cambridge University Press. Cambridge.
%
%\bibitem{Goetze} Götze, F. (1991). On the rate of convergence in the multivariate CLT. \textit{The
%Annals of Probability} \textbf{19}(2), 724-739.

\bibitem{ito}{ It{\^o}, K.} (1956).  {Spectral type of the shift transformation of differential processes  with
stationary increments}. {\it Trans. of the American Mathematical Society} {\bf 81},
252--263

\bibitem{Kab} {Kabanov, Y. (1975). }On extended stochastic integrals.\emph{\
}\textit{Theory of Probability and its Applications}{\ \textbf{20}},{\
710-722. }

\bibitem{KypBook} {Kyprianou, A.E. (2006).} \textit{Introductory Lectures
on Fluctuations of Lévy Processes with Applications}.
Springer-Verlag. Berlin Heidelberg New York.

\bibitem{lee1}{Lee, Y.-J. and Shih, H.-H. } (2004).
 The product formula of multiple L\'evy-It{\^o} integrals.
{\it Bull. Inst. Math. Acad. Sinica} {\bf 32}, 71--95.

\bibitem{LeSoUtVi} J. Leon, J.L. Solé, F. Utzet and J. Vives
(2002). On Lévy processes, Malliavin calculus and market models
with jumps. \textit{Finance and Stochastics} \textbf{6}, 197-225.

\bibitem{nieto05} {Nieto-Barajas, L.E. \textrm{and} Walker, S.G.}
(2005). A semi-parametric Bayesian analysis of survival data based on L\'{e}%
vy-driven processes. \textit{Lifetime Data Analysis} \textbf{11}, 529--543.

\bibitem{NouPe07} Nourdin, I.
and Peccati, G. (2008a). Stein's method on Wiener chaos. To appear
in: \textit{Probability Theory and Related Fields}.

\bibitem{NouPe08} Nourdin, I. and Peccati,
G. (2008b). Stein's method and exact Berry-Esséen bounds for
functionals of Gaussian fields. Preprint.

\bibitem{NouPeRev} Nourdin, I., Peccati, G. and Reveillac, A. (2008).
Multivariate normal approximation using Stein's method and Malliavin
calculus. Preprint.

\bibitem{NuaBook} Nualart, D.
(2006). \textit{The Malliavin calculus and related topics}. $2^{\rm
nd}$ Edition. Springer-Verlag. Berlin Heidelberg New York.

\bibitem{NO}
Nualart, D. and Ortiz-Latorre, S. (2007). \rm Central limit
theorems for multiple stochastic integrals and Malliavin calculus.
\textit{Stochastic Processes and their Applications} {\bf
118}(4), 614-628.

\bibitem{NP}
Nualart, D. and Peccati, G. (2005). Central limit theorems for
sequences of multiple stochastic integrals. \it The Annals of Probability {\bf
33} \rm (1), 177-193.

\bibitem{NualartVives}
{Nualart, D. \textrm{and} Vives, J.} (1990). Anticipative calculus
for the Poisson process based on the Fock space .
\textit{S\'{e}minaire de Probabilit\'{e}s XXIV}, LNM
\textbf{1426}, 154--165, Springer, Berlin.
%
%\bibitem{Ogura} Ogura, H. (1972). Orthogonal functionals of the Poisson
%process. \textit{IEEE Trans. Inform. Theory IT}, \textbf{18}, 473-481.

\bibitem{PePrU} Peccati, G. and Pr\"{u}nster, I. (2008). Linear and
quadratic functionals of random hazard rates: an asymptotic analysis. To
appear in: \textit{The Annals of Applied Probability.}

\bibitem{PeTaq0} Peccati, G. and Taqqu, M.S. (2007). Stable convergence of $%
L^{2}$ generalized stochastic integrals and the principle of conditioning.
\textit{The Electronic Journal of Probability}, \textbf{12}, 447-480, n. 15
(electronic).

\bibitem{PecTaqINT} Peccati, G.\ and Taqqu, M.S. (2007). Limit
theorems for multiple integrals. Preprint.

\bibitem{PecTaq07} Peccati, G.\ and Taqqu, M.S. (2008). Central
limit theorems for double Poisson integrals. To appear in:
\textit{Bernoulli}

\bibitem{PT} Peccati, G. and Tudor, C.A. (2004). Gaussian limits for
vector-valued multiple stochastic integrals. In: \textit{S\'{e}minaire de
Probabilit\'{e}s XXXVIII}, 247-262, Springer Verlag

\bibitem{Picard1996} Picard, J. (1996). Formules de dualit\'{e}s sur l'espace de Poisson. \textit{Ann. Inst. H. Poincar\'{e} B}, {\bf 32}(4), 509-548.

\bibitem{Reinert_sur} Reinert, G. (2005). Three general approaches to Stein's
method. In: \textit{An introduction to Stein's method}, 183-221.
Lect. Notes Ser. Inst. Math. Sci. Natl. Univ. Singap. \textbf{4},
Singapore Univ. Press, Singapore.

%\bibitem{RoWa} Rota, G.-C. and Wallstrom, C. (1997). Stochastic integrals: a
%combinatorial approach. \textit{The Annals of Probability} \textbf{25}(3),
%1257-1283.

\bibitem{Sato} Sato, K.-I. (1999). \textit{L\'{e}vy Processes and Infinitely
Divisible Distributions. }Cambridge Studies in Advanced Mathematics \textbf{%
68}. Cambridge University Press.

\bibitem{SolUtzViv} Solé, J.L., Utzet, F. and Vives, J. (2007). Canonical Lévy
processes and Malliavin calculus. \textit{Stochastic Processes
and their Applications} {\bf 117}, 165-187.

\bibitem{Stein_orig} Stein, Ch. (1972). A bound for the error in the normal approximation
to the distribution of a sum of dependent random variables. In:
\textit{Proceedings of the Sixth Berkeley Symposium on
Mathematical Statistics and Probability, Vol. II: Probability
theory}, 583-602. Univ. California Press, Berkeley, Calif..

\bibitem{Stein_book} Stein, Ch. (1986). \textit{Approximate computation of expectations.}
Institute of Mathematical Statistics Lecture Notes - Monograph
Series, \textbf{7}. Institute of Mathematical Statistics, Hayward,
CA.

\bibitem{Surg1984} Surgailis, D. (1984). On multiple Poisson stochastic
integrals and associated Markov semigroups. \textit{Probabilty and Mathematical Statistics}
\textbf{3}(2), 217-239.

\bibitem{Sur} Surgailis, D. (2000). CLTs for Polynomials of Linear Sequences:
Diagram Formulae with Applications. In: \textit{Long Range Dependence},
111-128, Birkh\"{a}user.

\bibitem{SUr2} Surgailis, D. (2000). Non-CLT's: U-Statistics, Multinomial
Formula and Approximations of Multiple Wiener-It\^{o} integrals. In: \textit{%
Long Range Dependence}, 129-142, Birkh\"{a}user.

%\bibitem{WolTaq} Wolpert, R.L. and Taqqu, M.S. (2005). Fractional
%Ornstein-Uhlenbeck L\'{e}vy Processes and the Telecom Process: Upstairs and
%Downstairs. \textit{Signal Processing }\textbf{85}(8), 1523-1545.
%
%\bibitem{WU} Wu, L. (1988). Construction de l'opérateur de Malliavin sur l'espace de Poisson.
%In: \textit{Séminaire de Probabilités XXI}, LNM \textbf{1247},
%100-113.

\bibitem{XUE} Xue, X.-H. (1991). On the principle of conditioning and convergence to mixtures of distributions for
sums of dependent random variables. \textit{Stochastic Processes and their Applications} \textbf{37}(2), 175-186.
\end{thebibliography}
\end{document}